\documentclass[11pt,a4paper,oneside]{amsart}
\usepackage{deleq}
\usepackage{cite}
\usepackage{url}
\usepackage{amssymb}
\usepackage{float}
\usepackage{a4}

\newcommand{\loc}{\textrm{loc}}
\newcommand{\mcH}{{\mycal H}}
\newcommand{\zmcH}{\;\mathring{\!\!\!\mycal H}}
\newcommand{\tg}{{\widetilde{g}}}
\newcommand{\can}{g_{\cerc^{n-1}}}%{\mathrm{can}}

\newcommand{\mzD}{D}
\newcommand{\pM}{\partial M }

\newcommand{\bel}[1]{\begin{equation}\label{#1}}
\newcommand{\dga}{\sigma_{g_a}}
\newcommand{\dgo}{\sigma_{g_1}}
\newcommand{\dgt}{\sigma_{g_2}}
\newcommand{\phot}{\phi_{12}}
\newcommand{\phto}{\phi_{21}}
\newcommand{\doto}[1]{\frac{d #1}{d s}}

\DeclareFontFamily{OT1}{rsfs}{}
\DeclareFontShape{OT1}{rsfs}{m}{n}{ <-7> rsfs5 <7-10> rsfs7 <10->
rsfs10}{} \DeclareMathAlphabet{\mycal}{OT1}{rsfs}{m}{n}
\swapnumbers

\theoremstyle{plain}
\newtheorem{theo}{Theorem}[section]
\newtheorem{adc}[theo]{Asymptotic decay conditions}
\newtheorem{lemm}[theo]{Lemma}
\newtheorem{pro}[theo]{Proposition}
\newtheorem{Proposition}[theo]{Proposition}

\theoremstyle{definition}

\theoremstyle{remark}
\newtheorem{remk}[theo]{Remark}

\numberwithin{equation}{section}

\newcommand{\bM}{{\overline{M}}}
\newcommand{\piM}{\partial_\infty M}
\newcommand{\piMx}{\partial_\infty \Mext}
\newcommand{\B}{{\mathbb B}}

\newcommand{\Z}{{\mathbb Z}}
\newcommand{\R}{{\mathbb R}}

\newcommand{\ra}{\rangle}
\newcommand{\la}{\langle}

\newcommand{\p}{\partial}

\newcommand{\br}{\bar{r}}
\newcommand{\bv}{{\bar{v}}{}}
\newcommand{\Eq}[1]{Equation~\eqref{#1}}
\newcommand{\Eqsone}[1]{Equations~\eqref{#1}}
\newcommand{\Eqs}[2]{Equations~\eqref{#1}-\eqref{#2}}

\newcommand{\tr}{\mathrm{tr}}
\newcommand{\cerc}{{\mathbb S}}

\newcommand{\Ric}{\mathrm{Ric}}
\renewcommand{\geq}{\geqslant}
\renewcommand{\leq}{\leqslant}
\newcommand{\dirac}{\mathfrak{D}}

\newcommand{\Int}{\mathrm{int}}
\newcommand{\Mext}{M_{{\text{\scriptsize\rm ext}}}}
\newcommand{\zh}{\breve{h}}

\newcommand{\zD}{{\mathring D}}
\newcommand{\Vect}{\mathrm{Vect}}
\newcommand{\nn}{\nonumber}
\newcommand{\tb}{{\tilde b}}

\newcommand{\cNb}{\mathcal{N}_b}
\newcommand{\cNbz}{\mathcal{N}_{b_0}}
\newcommand{\cNbone}{\mathcal{N}_{b_1}}
\newcommand{\cNbtwo}{\mathcal{N}_{b_2}}
\newcommand{\ourU}{\mathbb U}

\newcommand{\be}{\begin{equation}}
\newcommand{\ee}{\end{equation}}

\newcounter{mnotecount}[section]
\renewcommand{\themnotecount}{\thesection.\arabic{mnotecount}}
\newcommand{\mnote}[1]%{}
{\protect{\stepcounter{mnotecount}}$^{\mbox{\footnotesize  $%\!\!\!\!\!\!\,
      \bullet$\themnotecount}}$ \marginpar{\raggedright\tiny\em
    $\!\!\!\!\!\!\,\bullet$\themnotecount: #1} }

\newcommand{\eq}[1]{(\ref{#1})}

\begin{document}

\title[The mass of asymptotically hyperbolic manifolds]{The mass of
asymptotically hyperbolic Riemannian manifolds}
\author{Piotr T. Chru\'sciel}\
\address{D\'epartement de Math\'ematiques \\ UMR 6083 du CNRS\\
Universit\'e de Tours\\ Parc de Grandmont \\ F-37200 Tours, France}
\email{chrusciel@univ-tours.fr}
\author{Marc Herzlich}
\address{Institut de Math\'e\-matiques et Mod\'elisation de
Montpellier\\ UMR 5149 du CNRS\\ Universit\'e Montpellier~II\\
F-34095~Montpellier Cedex 5, France}
\email{herzlich@math.univ-montp2.fr}
\thanks{The first author is supported in part by the Polish Research
Council grant KBN 2 P03B 073 15. The second author is a member of the
EDGE Research Training Network HPRN-CT-2000-00101 of the European
Union and is also supported in part by the ACI
program of the French Ministry of Research.}

\begin{abstract}
We present a set of global invariants, called ``mass integrals",
which can be defined for a large class of asymptotically
hyperbolic Riemannian manifolds. When the ``boundary at infinity"
has spherical topology one single invariant is obtained, called
the mass; we show positivity thereof. We apply the definition to
conformally compactifiable manifolds, and show that the mass is
completion-independent. We also prove the result, closely related
to the problem at hand, that conformal completions of conformally
compactifiable manifolds are unique.
\end{abstract}


\maketitle

%\Input{todo2}
\section*{Introduction.}\label{Si} Let $(M,b)$ be a smooth
$n$-dimensional Riemannian manifold,
$n\ge 2$ and let $\cNb$ denote the set of functions $V$ on $M$
such that
\begin{eqnarray}
&  \label{eq:1}
  \Delta_b V + \lambda V =0\;,
& \\ & \label{eq:2} \zD_i\zD_j V = V( \Ric(b)_{ij} - \lambda
b_{ij})\;, &
\end{eqnarray}
for some constant $\lambda<0$. Here $\Ric(b)_{ij}$ denotes the
Ricci tensor of the metric $b$, $\zD$ the Levi-Civita connection
of $b$, and $\Delta_b:=b^{k\ell}\zD_k\zD_\ell $ is the Laplacian
of $b$. %We shall suppose that $\lambda<0$; s
Rescaling $b$ if necessary, we can without loss of generality
assume that
$$\lambda=-n \quad \textrm{so that} \quad R_b = b^{ij}\Ric(b)_{ij}
= -n(n-1)\;.$$
$(M,b)$ will be called {\em static\/} if
$$ \cNb \ne \{0\}\;.$$
This terminology is motivated by the fact that for every
$V\in\cNb$  the Lorentzian metrics defined on
$\R\times (M\setminus\{V=0\})$ by the formula
\begin{equation}
  \label{eq:3}
\gamma = - V^2 dt^2 + b
\end{equation}
are static solution of the Einstein equations, $\Ric({\gamma})
=\lambda\gamma$
(and the Riemannian metrics $V^2 dt^2 + b$ are actually Einstein as well).

The object of this work is to present a set of global invariants,
constructed using $\cNb$, for metrics which are asymptotic to a
class of static metrics. The model case of interest is the
hyperbolic metric, which is static in our sense: the set $\cNb$ is
then linearly isomorphic to $\R^{n+1}$; however, other classes of
metrics will also be allowed in our framework. The invariants
introduced here stem from a Hamiltonian analysis of general
relativity, and part of the work here is a transcription to a
Riemannian setting of the Lorentzian analysis in
\cite{ChNagyATMP}. Related definitions have been given recently by
Wang~\cite{Wang} (with a spherical conformal infinity) and
Zhang~\cite{Zhang:hpet} (in dimension $3$, again with spherical
asymptotic geometry), under considerably more restrictive
asymptotic and global conditions.  Wang's definition of mass for
asymptotically hyperbolic manifolds~\cite{Wang} coincides with
ours, when his much stronger asymptotic decay conditions are
satisfied (note, however, that his proof  of the geometric
character of the mass is incomplete, as he ignores the possibility
of existence of inequivalent conformal completions). Moreover, the
hypothesis of \cite{Wang,Zhang:hpet} that $(M,g)$ has compact
interior is replaced by that of completeness; this strengthening
of the positivity theorem is essential when the associated
Lorentzian space-time contains ``degenerate'' event horizons.

 This work is organised as follows. In
Section~\ref{S2} we review a few static metrics, and discuss their
properties relevant to the work here. In Section~\ref{S3} we
define the ``mass integrals'', we make precise the classes of
metrics considered, and we show how to obtain global invariants
out of the mass integrals. We show that our boundary conditions
are sharp, in the sense that their weakening leads to mass
integrals which do not provide  geometric invariants. It is
conceivable that the strengthening of some of our conditions could
allow the weakening of some other ones, leading to geometric
invariants for other classes of manifolds; we expect that such a
mechanism occurs for the Trautman-Bondi mass of asymptotically
hyperboloidal manifolds. In several cases of interest one obtains
a single invariant, which we call the mass of $(M,g)$, but more
invariants are possible depending upon the topology of the
``boundary at infinity" $\piM$ of $M$
---this is determined by the number of invariants of the action
of the group of isometries of $b$ on $\cNb$, see
Section~\ref{Sminv} for details. Regardless of this issue, we
emphasise that we only consider the ``global charges'' of
\cite{ChNagyATMP} related to Killing vectors which are normal to
the level sets of $t$ in  the space-time metric \eq{eq:3}: the
remaining space-time invariants of \cite{ChNagyATMP}, associated,
{\em e.g.},\/ to ``rotations'' of $M$, involve the extrinsic
curvature of the initial data hypersurface and are of no concern
in the purely Riemannian setting here. In Section~\ref{S4} we
prove positivity of the mass so obtained for metrics asymptotic to
the hyperbolic one. As a corollary of the positivity results we
obtain a new uniqueness result for anti-de Sitter space-time,
Theorem~\ref{Ctheo}. We also consider there the case of manifolds
with a compact inner boundary. In Section~\ref{Sccm} we show how
to define mass for a class of conformally compactifiable
manifolds. The question of geometric invariance of the mass is
then closely related to the question of uniqueness of conformal
completions; in Section~\ref{S5} we prove that such completions
are indeed unique.

\bigskip

\section{The reference metrics}
\label{S2}
 Throughout this paper we will assume that the manifold
contains a region $\Mext\subset M$ together with a  diffeomorphism
\be\label{Nman}\Phi^{-1}:\Mext \to [R,\infty)\times N\;,\ee
where $N$ is a compact boundaryless manifold, such that the reference
metric $b$ on $\Mext$ takes the form
\be
\label{cm1} \Phi^*b=\frac{dr^2}{r^2+k}+r^2\zh=: b_0\;, \ee with
$\zh$
---a Riemannian metric on $N$ with scalar  curvature $R_{\zh}$ and the
constant $k$ equal to
\be\label{cm1.1}
R_{\zh}=
\begin{cases}(n-1)(n-2)k\;, \quad k\in\{0,\pm 1\}
\;, & \mbox{if }\ n>2\;, \cr 0\;, \quad k=1 \;, & \mbox{if }\
n=2\;, \end{cases}\ee (recall that the dimension of $N$ is
$(n-1)$); here $r$ is a coordinate running along the $[R,\infty)$
factor of $[R,\infty)\times N$. There is some freedom in the
choice of $k$ in \eq{cm1} when $n=2$, associated with the range of
the angular variable $\varphi$ on $N=S^1$ (see the discussion in
Remark~\ref{Rdim2} below) and we make the choice $k=1$, as it
corresponds to the usual form of the two-dimensional hyperbolic
space.

When $(N,\zh)$  is the unit round $(n-1)$--dimensional sphere
$(\cerc^{n-1},\can)$, then $b$ is the hyperbolic metric.
Equations~\eq{cm1} and \eq{cm1.1} imply that the scalar curvature
$R_b$ of the metric $b$ is constant:
 $$ R_b=-
n(n-1)\;.$$ Moreover the metric $b$ will be Einstein if and only
if $\zh$ is. We emphasise that for all our purposes we only need
$b$ on $\Mext$, and we continue $b$ in an arbitrary way to
$M\setminus\Mext$ whenever required.

The cases of main interest seem to be those where $\zh$ is a space form
--- it then follows from the results in \cite[Appendix~B]{ChNagyATMP}
that we have
\begin{eqnarray}&
\label{cm2} k=0,-1 \quad \Longrightarrow \quad \cNbz=
\Vect\{\sqrt{r^2+k}\}\;,&
\\%&&
\label{cm3.1}%\\\nn
&
 k=1\;,\ (N,\zh)=(\cerc^{n-1},\can)  \quad \Longrightarrow \quad \cNbz =
\Vect\{V_{(\mu)}\}_{\mu=0,\ldots,n}\;,  &\\&
V_{(0)}=\sqrt{r^2+k}\;, \quad V_{(i)}=x^i\;,&\label{cm3.2}
\end{eqnarray}
with the usual identification of $[R,\infty)\times \cerc^{n-1}$
with a subset of $\R^n$ in \eq{cm3.2}. However, we shall not
assume that $\zh$ is a space-form, or that \Eqs{cm2}{cm3.2} hold
unless explicitly stated.

 For the purposes of
Section~\ref{Sccm} we note the following conformal representation
of the metrics \eq{cm1}: one replaces the coordinate $r$ by a
coordinate $x$ defined as \be\label{radcoord} x = \frac 2
{r+\sqrt{r^2+k}} \quad \Longleftrightarrow \quad r =
\frac{1-kx^2/4}{x}\;,\ee which brings $b$ into the form
\begin{eqnarray}
\label{cm5} b & = & x^{-2}\left(dx^2+(1-kx^2/4)^2\zh\right)
\\
\nn &=:& x^{-2}\tb\;,
\end{eqnarray}
with $\tb$ -- a metric smooth up to boundary on
$\{x\in[0,x_R]\times N\}$, for a suitable $x_R$.

\section{The mass integrals} \label{S3} Let $g$ and $b$ be two
Riemannian metrics on a manifold $M$, and let $V$ be any function
there. We set
\begin{equation}
  \label{eq:3.1}
 e_{ij} :=g_{ij}-b_{ij}
\end{equation}
(the reader is warned that the tensor field $e$ here is \emph{not}
a direct Riemannian counterpart of the one in~\cite{ChNagyATMP};
the latter makes appeal to the \emph{contravariant} and not the
\emph{covariant} representation of the metric tensor). As before
we denote by $\zD$ the Levi-Civita connection of $b$, and we use
the symbol $R_f$ to denote the scalar curvature of any metric $f$.
The basic identity from which our mass integrals arise is the
following:
\begin{eqnarray}
  \label{eq:3.2} &
  \sqrt{\det g} \;V (R_g-R_b)  =  \partial_i \left(\ourU^i(V)\right) + \sqrt{\det g}
\;(  \rho + Q)\;, &
\end{eqnarray}
where
\begin{eqnarray} \label{eq:3.3} & {}\ourU^i (V):=  2\sqrt{\det
g}\;\left(Vg^{i[k} g^{j]l} \zD_j g_{kl}
%\phantom{xxxxxxxx} &\\&\nn \phantom{xxxxx\ourU^i :=2}
%+
%(n-1)\zD^{i}V(\sqrt{\det g} -\sqrt{\det b})\;,&
+D^{[i}V %(\sqrt{\det g}
g^{j]k} e_{jk}\right)
% -\sqrt{\det b}b^{j]k}
%b_{jk} )
\;,&
\\
\label{eq:3.4} & \rho := (-V\Ric(b)_{ij} +\zD_i \zD_j V -\Delta_b
V b_{ij}) g^{ik} g^{j\ell}
  e_{k\ell}
\;,&
\\
\label{eq:3.5} & Q:= V(g^{ij} - b^{ij} + g^{ik}g^{j\ell}
e_{k\ell})\Ric(b)_{ij} +Q'\;.
\end{eqnarray} Brackets over a symbol denote
anti-symmetrisation, with an appropriate numerical factor ($1/2$
in the case of two indices). Here $ Q'$ denotes an expression
which is bilinear in $e_{ij}$ and $\zD_k e_{ij}$, linear in $V$,
$dV$ and Hess$V$, with coefficients which are constants in an ON
frame for $b$. The idea behind this calculation is to collect all
terms in $R_g$ that contain second derivatives of the metric in
$\partial_i \ourU^i$; in what remains one collects in $\rho$ the
terms which are linear in $e_{ij}$, while the remaining terms are
collected in $Q$; one should note that the first term at the
right-hand-side of \eq{eq:3.5} does indeed not contain any terms
linear in $e_{ij}$ when Taylor expanded at $g_{ij}=b_{ij}$. The
mass integrals will be flux integrals
--- understood as a limiting process --- over the ``boundary at
infinity of $M$" of the vector density $\ourU^i$

In general relativity a normalising factor $1/16\pi$, arising from
physical considerations, is usually thrown in into the definition
of $\ourU^i$. From a geometric point of view this seems purposeful
when the boundary at infinity is a round two dimensional sphere;
however, for other topologies and dimensions,  this choice of
factor  does not seem very useful, and for this reason we do not
include it in the definition.

We note that the linearisation of the mass integrands $\ourU^i$
coincides with the linearisation of the charge integrands
of~\cite{ChNagyATMP} evaluated for the Lorentzian metrics
$^4b=-V^2dt^2+b$, $^4g=-V^2dt^2+g$, with $X=\partial_t$, on the
hypersurface $t=0$; however the integrands do not seem to be
identical. Nevertheless, under the conditions of
Theorem~\ref{Tinv} the resulting numbers coincide, because under
the asymptotic conditions of Theorem~\ref{Tinv} only the
linearised terms matter.

The convergence of the mass integrals requires appropriate
boundary conditions, which are defined using the following
orthonormal frame $\{f_i\}_{i=1,n}$ on $\Mext$:
\be\label{m1} \Phi^{-1}_*f_i = r^{-1}\epsilon_i\;, \quad
i=1,\ldots,n-1\;,\quad\Phi^{-1}_*f_n = \sqrt{r^2+k}
\;\partial_r\;, \ee where the $\epsilon_i$'s form an orthonormal
frame for the metric $\zh$. We moreover set \be \label{m2}
g_{ij}:=g(f_i,f_j)\; .\ee
\begin{adc} We shall require:
\begin{deqarr}\arrlabel{m3}\label{m3a}
& \int_{\Mext} \left( \sum_{i,j} |g_{ij}-\delta_{ij}|^2 + \sum_{i,
j,k} |f_k(g_{ij})|^2 \right)r\circ\Phi\;d\mu_g<\infty\;,&\\
& \int_{\Mext}
|R_g-R_b|\;r\circ\Phi\;d\mu_g<\infty\;,\label{m3b}&\end{deqarr}
\begin{equation}
\label{m0} \exists \ C > 0 \ \textrm{ such that }\
C^{-1}b(X,X)\le g(X,X)\le Cb(X,X)\;. \end{equation}
\end{adc}

\medskip

For the $V$'s of \Eqsone{cm2} or \eq{cm3.2} we have
\begin{equation}
\label{Vcondition} V=O(r)\;, \quad \sqrt{b^\#(dV,dV)}=O(r)\;,
\end{equation}
where $b^\#$ is the metric on $T^*M$ associated to $b$, and this
behavior will be assumed in what follows:

\begin{pro}
\label{P1} Let the reference metric $b$ on $\Mext$ be of the form
\eqref{cm1}, suppose that $V$ satisfies \eq{Vcondition}, and
assume that $\Phi$ is such that Equations \eq{m3}-\eq{m0} hold.
Then for all $V\in\cNbz$ the limits \be \label{mi}
H_\Phi(V):=\lim_{R\to\infty} \int_{r=R} \ourU^i(V\circ \Phi^{-1})
dS_i \ee exist, and are finite.
\end{pro}

The integrals~\eq{mi} will  be referred to as \emph{the mass
integrals}.

\medskip

\proof For any $R_1,R_2$ we have \be \int_{r=R_1} \ourU^i
dS_i = \int_{r=R_2} \ourU^i dS_i + \int_{[R_1,R_2]\times N}
\partial_i \ourU^i \,d^n x\;, \ee and the result follows from
\eqref{eq:3.2}-\eqref{eq:3.5}. \qed

Under the conditions of Proposition~\ref{P1}, the
integrals~\eqref{mi} define a linear map from $\cNbz$ to $\R$.
Now, each map $\Phi$ used in \eq{cm1} defines in general a
different background metric $b$ on $\Mext$, so that the maps
$H_\Phi$ are potentially dependent upon  $\Phi$.
(It should be clear that, given a
fixed $\zh$, \eq{mi} does not depend upon the choice of the frame
$\epsilon_i$ in \eq{m1}.) It turns out that this dependence can be
controlled, as follows:

\begin{theo}\label{Tinv}
Consider two maps $\Phi_a$, $a=1,2$,
satisfying \eq{m3} together with \be \label{m5}\sum_{i,j}
|g_{ij}-\delta_{ij}| + \sum_{i, j,k} |f_k(g_{ij})|= \begin{cases}
o(r^{-n/2})\;, & \mbox{if }\ n>2\;, \cr O(r^{-1-\epsilon})\;, &
\mbox{if }\ n=2\;, \mbox{for some $\epsilon>0$}.\end{cases} \ee
Then there exists an isometry $A$ of $b_0$, defined perhaps only
for $r$ large enough, such that \bel{Htransf} H_{\Phi_2}(V)=
H_{\Phi_1}(V\circ A^{-1}) \ee
\end{theo}

\proof The arguments of the proof of Theorem~\ref{Tinv} follow
closely those given at the beginning of Section~4 and in Section~2
of~\cite{ChNagyATMP}, we need, however, to adapt some of the
necessary ingredients to our different setup here. The conclusion
of Proposition~\ref{pscc} below, which holds for all manifolds
$(N,\zh)$ considered here, enables us to use Theorem 3.3 (2) of
\cite{ChNagyATMP}: if $\Phi_1$ and $\Phi_2$ are two maps as above
satisfying the decay assumptions (\ref{m3})-(\ref{m0}) and
(\ref{m5}) with respect to (isometric) reference metrics $b_1:=
(\Phi^{-1}_1)^*b_0$ and $b_2:= (\Phi_2^{-1})^*b_0$, then there
exists an isometry $A$ of the background metric $b_0$, defined
perhaps only for $r$ large enough, as made clear in
Proposition~\ref{pscc}, such that $$\Phi_2 - \Phi_1\circ A =
o(r^{-n/2})\;$$ One also has a similar
--- when appropriately formulated in terms of $b$-ortho\-nor\-mal
frames, as in \cite{ChNagyATMP}
--- decay of first two derivatives.

It follows directly from the definition of $H_{\Phi}$ that $$
H_{\Phi_1\circ A} (V)= H_{\Phi_1}(V\circ A^{-1})\;.$$ In order to
establish \eq{Htransf} it remains to show that \bel{Htransf2}
H_{\Phi_1\circ A} (V)= H_{\Phi_2}(V)\;. \ee Now, Corollary~3.5 of
\cite{ChNagyATMP} shows that $\Phi_1\circ A$ has the same decay
properties as $\Phi_1$, so that --- replacing $\Phi_1$ by
$\Phi_1\circ A$ --- to prove  \eq{Htransf2} it remains to consider
two maps $\Phi_1^{-1} =(r_1,v^A_1)$ and $\Phi_2^{-1}=(r_2,v^A_2)$
(where $v^A$ denote abstract local coordinates on $N$) satisfying
\begin{equation}
\begin{split}
r_2 & = r_1 + o(r_1^{1-\frac{n}{2}}) ,\\
v^A_2 & = v^A_1 + o(r_1^{-(1+\frac{n}{2})}) ,
\end{split}\end{equation}
together with elements $V_1:=V\circ \Phi_1^{-1}$ of $\cNbone$ and
$V_2:=V\circ \Phi_2^{-1}$ of $\cNbtwo$ having the {\it same
expression} in the first or the second system of coordinates.
Local coordinates $v^A$ might not be defined on the whole of $M$;
we shall remove this problem by embedding the manifold $N$ in
$\R^{2(n-1)}$, so that local coordinates are turned into global
coordinates. This has no effect in the sequel of the proof but
enables us to consider a well-defined vector field
\[ \zeta = (r_2 - r_1)\frac{\partial}{\partial r_1} + \sum_A (v^A_2 - v^A_1)
\frac{\partial}{\partial v^A_1} \ , \] defined only along $M$, and
tangent to $M$. The decay estimates above imply that $\zeta =
o(r^{-n/2})$ in the reference metric $b_1$; by Theorem 3.3 (2) of
\cite{ChNagyATMP} the same holds for its first two
$\zD$-derivatives. Elementary calculations show then that
\begin{equation}\label{eq:newb}
b_2 = b_1 + \mathcal{L}_{\zeta}b_1 + o(r^{-n})\ , \ \ V_2 = V_1 +
\zD_iV\zeta^i  + o(r^{1-n})\;,
\end{equation}
together with their first derivatives. Hence, to leading order in
powers of $r\approx r_1$, everything behaves as if we were
considering a first order variation of metrics through the action
of the flow of the vector field $\zeta$.

We shall now show that $H_{\Phi_1}(V) = H_{\Phi_2}(V)$.
For the purpose of the
calculations that follow it will be easier to replace the local
integrand $\ourU$ by the following one,
\begin{equation}\label{eq:newU}
\ourU^i = \sqrt{\det\, b}\left( - V \zD_jg^{ij} + V b^{ij}b_{kl}
\zD_jg^{kl} + 2\,\zD^{[i}Vb^{j]k}e_{jk} \right)\;,
\end{equation}
which yields the same limit at infinity when integrated on an
element $V$ of $\cNb$ on larger and larger spheres (strictly
speaking, we should not denote them by the same letter $\ourU$,
since they are different vector densities which give identical
results only after an integration process; we shall however do so
since expression (\ref{eq:newU}) will only be used in the course
of the current proof; we emphasise that the definition
(\ref{eq:3.3}) is used in all other places in the paper).

We now compute the variation of $\ourU$ when passing from the
asymptotic map $\Phi_1$ (with reference metric $b_1$ and function
$V_1$) to the second map $\Phi_2$ (with reference metric $b_2$ and
function $V_2$). {}From Equation (\ref{eq:newb}), we deduce
\begin{equation}
\ourU^i_2 - \ourU^i_1 = \delta\ourU^i + o(r^{1-n}),
\end{equation}
where $\delta\ourU^i$ is obtained by linearisation in $\zeta$ at
$g=b$ and will be computed below, while the remainder terms decay
sufficiently fast so that they do not contribute when integrated
at infinity against either $b_1$ or $b_2$. It remains to show that
$\delta\ourU^i$ does not contribute either when integrated at
infinity. In Equation (\ref{eq:newU}), the only terms that
contribute {\it a priori} to the variation of $\ourU$ are the
following: $b^{ij}$, $b_{kl}$, $\sqrt{\det\, b}$, $\zD$ and $V$,
but the decay estimates (\ref{eq:newb}) show at first glance that
only the variation of $\zD$ will contribute to the first-order
term $\delta\ourU$. We now compute it using Formulas 1.174
of~\cite{Besse}. In all what follows, we denote $b=b_1$ and
$V=V_1$. Then,
\begin{equation}\begin{split}
\delta\ourU^i = & \sqrt{\det\, b}\left( - V \zD_k\zD^k \zeta^i  + V \zD_k
\zD^i\zeta^k - 2V \Ric(b)^i{}_k\zeta^k \right) \\
& \ + \sqrt{\det\, b}\left( - 2 (\zD^iV)\zD_k\zeta^k + (\zD_kV)
\zD^i\zeta^k + (\zD^kV)\zD_k\zeta^i \right)\;.\end{split}
\end{equation}
Fortunately, this will appear to be the sum of a divergence term
plus lower order terms. The first step is to use the following elementary
facts:
\begin{equation}\begin{split}
- V \zD_k\zD^k\zeta^i & = - \zD_k(V\zD^k\zeta^i) + (\zD_kV)\zD^k\zeta^i,\\
V \zD_k\zD^i\zeta^k & = \zD_k(V\zD^i\zeta^k) - (\zD_kV)\zD^i\zeta^k,
\end{split}\end{equation}
which yield
\begin{equation}\label{twosp}\begin{split}
\delta\ourU^i \ = & \ 2\sqrt{\det\, b}\left( (\zD_kV)\zD^k\zeta^i
- (\zD^iV)\zD_k\zeta^k - V \Ric(b)^i{}_k\zeta^k \right) \\
& \ + \textrm{ divergence term}\; .\end{split}
\end{equation}
Each of the first two terms in the right-hand side may be transformed
with
\begin{equation}\begin{split}
(\zD_kV)\zD^k\zeta^i & = \zD^k(\zeta^i\zD^k V) - (\zD^k\zD_kV)\zeta^i, \\
- (\zD^iV)\zD_k\zeta^k & = - \zD_k(\zeta^k\zD^iV) + (\zD_k\zD^iV)\zeta^k,
\end{split}\end{equation}
and one may also use that $V$ is an element of $\cNb$ to conclude
that
\begin{equation}
\ourU^i_2 - \ourU^i_1 =  \textrm{ divergence term} + o(r^{1-n}).
\end{equation}
This establishes the covariance of the mass  functional. \qed

\medskip

\begin{remk} For the purpose of explicit calculations we
 note that under~\eq{m5} the mass
integral $H_{\Phi}(V_{(0)})=H_{\Phi}(\sqrt{r^2+k})$ can be written
as:
 \begin{eqnarray}
   \label{massequation1}
 && \\ \nn
   \displaystyle H_{\Phi}(V_{(0)})
&=&\lim_{R\to\infty}  (R^2+k)\times
\\ \nn &&
   \displaystyle \int_{\{r=R\}} \left(-\sum_{i=1}^{n-1}\left\{\frac {
          \partial e_{ii}}{\partial r}+ \frac {
          k e_{ii}}{ r(r^2+k)}\right\}+\frac {(n-1)e_{nn}}{r}%+e^{\hA\hA}
   \right) d^{n-1}\mu_{h}\;,
 \end{eqnarray}
assuming that the right-hand-side of \eq{massequation1} converges.
Here $d^{n-1}\mu_{h}$ is the Riemannian measure associated with
the metric $h$ induced on the level sets of the function $r$.
\end{remk}
%Further,
%$g^{ij}:=g^{\#}(\theta^i,\theta^j)$, where $g^{\#}$ is the metric
%on $T^*M$ associated to $g$, while the $\theta^i$'s form a
%co-frame dual to the $f_i$'s defined in \Eq{m1}.

\medskip

\begin{remk}
The conditions \eq{m5} are sharp, in the following sense: let $g$
be the standard hyperbolic metric, thus in a coordinate system
$(\br,\bv^A)$, where the $\bv^A$'s are local coordinates on
$\cerc^{n-1}$, we have \be \label{cm1bar}
g=\frac{d\br^2}{\br^2+k}+\br^2\zh\;.\ee  Let, for sufficiently
large $r$,
 $\Phi_\gamma^{-1}(r,v^A)=(\br(r,v^A),\bv^B(r,v^A))$ be given by
the formula \be \label{m6} \br=r+\gamma r^{1-n/2}\;,\qquad
\bv^A=v^A\;, \ee where $\gamma$ is a constant. Then
$H_{\Phi_\gamma}(\sqrt{r^2+k})$ \emph{does} depend upon
$\gamma$: in order to see that, consider any
background metric of the form
$$ b= a^2(r)dr^2 + r^2 \zh\;,$$
and let $g$ satisfy \bel{ansatz} g= g_{nn}(r,v^A)a^2(r)dr^2+
c(r,v^A)r^2 \zh\;,\ee for some differentiable functions $g_{nn}$
and $c$. One finds
\begin{eqnarray} \label{mansatz} &&\\\nn \ourU^i
dS_i|_{r=\textrm{const}}&:=& \ourU^i
\partial_i \rfloor dr\wedge dv^1\wedge \cdots \wedge
dv^{n-1}|_{r=\textrm{const}}
\\ \nn &=& \displaystyle
\frac{(n-1)c^{(n-3)/2} r^{n-2}}{a(r)\sqrt{g_{nn}}} \left\{
V\left(g_{nn}-1 - r \frac
{\partial c}{\partial r}\right) %+\left(V+\phantom{r\frac{\partial V}{\partial v}}\right.
\right.\\
\nn && %\left.
\left.
 +\left(%\phantom{r\frac{\partial V}{\partial v}}\right.
 r\frac{\partial V}{\partial
r}-V\right)(c-1)\right\}\; \sqrt{\det \zh_{AB}}\;dv^1\wedge \cdots
\wedge dv^{n-1}\;.
\end{eqnarray}
Applying this formula to the above $g$ and $b$ one obtains
\bel{mnotinv}H_{\Phi_\gamma}(\sqrt{r^2+k})= \frac 14
(n+8)n(n-1)\gamma^2\mbox{\rm Vol}_{\can}(\cerc^{n-1})\;.\ee One
can also check that the numerical value of the linearised
expression \eq{massequation1} reproduces the right-hand-side of
\eq{mnotinv} for the metrics at hand, thus is again \emph{not}
invariant under \eq{m6}.
\end{remk}

\section{The mass}
\label{Sminv}  In the asymptotically flat case the mass is a
single number which one assigns to each end of $M$
\cite{Bartnik,ChErice}; it is then natural to enquire whether
there are some geometrically defined {\sl numbers} one can extract
out of the family of maps $H_{\Phi}$. This will depend upon
the structure of $\cNbz$ and we shall give here a few important
examples. Throughout this section we assume that $\cNbz\ne
\emptyset$.

{\flushleft\bf A}. The simplest case is that of the manifold
$(N,\zh)$ of \eq{Nman}-\eq{cm1} having a strictly negative Ricci
tensor,  with scalar curvature $R_{\zh}$ equal to
$-(n-1)(n-2)$, %is a space-form with sectional curvature $-1$,
so that $n\ge 3$ and $k=-1$ in \eq{cm1}. Similarly to the space
forms discussed in Section~\ref{S2}, $\cNbz$ is
then~\cite[Appendix~B]{ChNagyATMP} one dimensional:
\be\label{sm1} V\in \cNbz \ \Longleftrightarrow V=\lambda
V_{(0)}\;,\ \lambda\in \R\;, \quad V_{(0)}:=\sqrt{r^2+k} \;. \ee
The coordinate system of \eq{cm1} is uniquely defined, so is the
function $V_{(0)}$; the number \be\label{sm2} m:=
H_\Phi(V_{(0)})\;,\ee calculated using any
 $\Phi$ satisfying the conditions of Theorem~\ref{Tinv}, provides
the desired, $\Phi$-independent definition of \emph{mass} relative
to $b_0$, whenever \eq{m1} and \eq{m5} hold.

{\flushleft\bf B}. Consider, next, the case of a flat $(N,\zh)$
with $n\ge 3$, so that $k=0$ in \eq{cm1}. \Eq{sm1} holds again;
however, the coordinate $r$ is {\em not} anymore uniquely defined
by $b$, since \eq{cm1} is invariant under the rescalings
$$r\to ar\;,\quad \zh\to a^{-2} \zh\;, \quad a\in \R^*\;.$$
 This freedom can be gotten rid of by requiring, \emph{e.g.},
 $$ \mathrm{Vol}_{\zh}(N)=1\;;$$
 the number $m$ obtained then from \eq{sm2}, with $V_{(0)}$ as in
\eq{sm1}, provides the desired invariant.

{\flushleft\bf C}. The case $k=+1$ requires more work. Consider
first the case where $(N,\zh)=(\cerc^{n-1},\can)$, so that the
reference metric is the hyperbolic metric; it is convenient to
start with a discussion of $\cNb$ in two models of the hyperbolic
space: in the ball model, we consider the ball $\B = \{ x \in
\R^n, |x| <1\}$ endowed with the metric $b = \omega^{-2} \delta$,
where \be\label{omf} \omega = \frac{1}{2}(1 - |x|^2)\;,\ee and
$\delta$ is the flat Euclidean metric. From \eq{cm3.1}-\eq{cm3.2}
one finds that the set $\cNb$ defined in \eq{eq:1}-\eq{eq:2} is
the $(n+1)$-dimensional vector space spanned by the following
basis of functions
\begin{equation}\label{functions}
V_{(0)} = \frac{1+|x|^2}{1-|x|^2}\;,\qquad V_{(i)} =
\frac{2x^i}{1-|x|^2}\;,
\end{equation}
where $x^i$ is any of the Cartesian coordinates on the flat ball.
In geodesic coordinates around an arbitrary point in the
hyperbolic space, the hyperbolic metric is $b=dr^2 + \sinh^2(r)
g_{\cerc^{n-1}}$ and the above orthonormal basis of $\cNb$ may be
rewritten as $V_{(0)} = \cosh(r)$ and $V_{(i)} = \sinh(r)n^i$,
where $n^i$ is the restriction of $x^i$ to the unit sphere centred
at the pole.

The space $\cNb$ is naturally endowed with a Minkowski metric
$\eta$, with signature $(+,-,\cdots,-)$, issued from the action of
the group of isometries $O^+(n,1)$ of the hyperbolic metric ({\em
cf., e.g.\/}~\cite[Appendix~B]{ChNagyATMP}). The basis given above
is then orthonormal with respect to this metric, with the vector
$V_{(0)}$ being timelike, \emph{i.e.},\/
$\eta(V_{(0)},V_{(0)})>0$. We define the time orientation of
$\cNb$ using this basis --- by definition a timelike vector
$X^{(\mu)}V_{(\mu)}$ is future directed if $X^{(0)}>0$, similarly
for covectors; this finds its roots in a Hamiltonian analysis in
the associated Lorentzian space-time. Assuming that there exists a
map $\Phi$ for which the convergence conditions of
Proposition~\ref{P1} are satisfied, we set \be\label{pmu}
p_{(\mu)}:=H_\Phi(V_{(\mu)})\;.\ee Under isometries of $b$ the
$V_{(\mu)}$'s are reshuffled amongst each other under the usual
covariant version of the defining representation of the Lorentz
group $O^+(n,1)$. It follows that the number
\be\label{masssphere} m^2:
=|(p_{(0)})^2-\sum_{i=1}^n(p_{(i)})^2|\ee is a geometric
invariant, which provides the desired notion of mass for a
spherical asymptotic geometry. The nature of the action of
$O^+(n,1)$ on $\cNb$ shows that the only invariants which can be
extracted out of the $H_\Phi$'s are $m^2$ together with the causal
character of $p_{(\mu)}$ and its future/past pointing nature if
relevant. Under natural geometric conditions $p_{(\mu)}$ is
timelike future pointing or vanishing, see Section~\ref{S4} below.
For timelike $p_{(\mu)}$'s it appears natural to choose the sign
of $m$ to coincide with that of $p_{(0)}$, and this is the choice
we shall make.

Suppose, finally, that the manifold $(N,\zh)$  is the
quotient of the  unit round sphere $(\cerc^{n-1},\can)$ by a
subgroup $\Gamma$ of its group of isometries. For generic
$\Gamma$'s one expects the conformal isometry group of $(N,\zh)$
to be trivial, in which case all the integrals $p_{(\mu)}$ defined
by \Eq{pmu} define invariants. In any case, for non-trivial
$\Gamma$'s conformal isometries of $(N,\zh)$ are
isometries (for compact Einstein manifolds which are not
round spheres the group of conformal isometries coincides with the
group of isometries; this follows immediately, {\em e.g.,}\/ from
what is said in~\cite{Kuehnel}; P.T.C. is grateful to A.~Zeghib
and C.~Frances for useful comments concerning the structure of the
group of conformal isometries of quotients of spheres), and, in
addition to $m$, $p_{(0)}$ becomes then a geometric invariant.
Further invariants may occur depending upon the details of the
action of the group of isometries of $(N,\zh)$ on $\cNbz$.

\begin{remk}\label{Rdim2}
Our results also apply in dimension $n=2$. This might seem
somewhat surprising at first sight, because there is no direct
useful equivalent of asymptotic flatness and of the associated
notion of mass given by an ADM-type integral in dimension $2$:
when the scalar curvature is in $L^1$, the appropriate analogue of
mass is the deficit angle, as made precise by the Shiohama
theorem~\cite{Shiohama}. For
the metrics considered here the Shiohama theorem does not apply;
however, the metrics we study can be thought of as having a minus
infinite deficit angle, consistently with a naively understood
version of the Shiohama theorem --- the ratio of the length of
distance circles to the distance  from any compact set tends to
infinity as the distance does. Examples of metrics on $\Mext$
which satisfy our asymptotic conditions and have a well defined
non-trivial mass
--- with respect to a background given by  the $2$-dimensional
hyperbolic metric --- are provided by the Riemannian counterpart
of the generalised $(2+1)$-dimensional Kottler metrics,
\be\label{twoK} g=\frac{dr^2}{r^2-\eta} +r^2 d\varphi^2\;, \quad
\varphi\in[0,2\pi]\ \mbox{mod\ } 2\pi \;,\ee for some constant
$\eta\in\R$; $\eta=-1$ corresponds to the standard hyperbolic metric
(as pointed out by Ba\~nados, Teitelboim and Zanelli~\cite{BTZ}, for
positive $\eta$ the associated static Lorentzian space-times with
$V=\sqrt{r^2-\eta}$ can be extended to space-times containing a
black-hole region). The metrics~\eq{twoK} have
constant Gauss curvature equal to minus one for all $\eta\in\R$,
so the integral condition on $R_g$ in \eq{m3} holds with $$b
=\frac{dr^2}{r^2+1} +r^2 d\varphi^2\;;$$ the remaining conditions
arising from \eq{m3}, as well as \eq{m0} and \eq{m5}, are easily
checked. Applying formula~\eq{massequation1} one obtains
$$m= p_{(0)}=2\pi(1+\eta)$$
(the remaining $p_{(\mu)}$'s are zero by symmetry
considerations). For strictly negative $\eta$ there is a
sense in which $m$ is related to a deficit angle, as follows: a
coordinate transformation $r\to\lambda r$, $\varphi \to
\varphi/\lambda$, with $\lambda^2=-\eta$ brings the metric
\eq{twoK} to the standard hyperbolic space form, \be\label{twoK2}
g=\frac{dr^2}{r^2+1} +r^2 d\varphi^2\;, \quad
\varphi\in[0,2\pi/\lambda]\ \mbox{mod\ } 2\pi/\lambda \;,\ee
except for the changed range of variation of the angular variable
$\varphi$; that range will coincide with the standard one if and
only if the mass vanishes. For metrics asymptotic to \eq{twoK}
%with leading order
%behavior (\ref{twoK}), in the sense of \eq{m5} with $\epsilon>1$,
the geometric invariance of the mass should follow directly from
this deficit angle character.
% \ptc{argument added} follows directly
%from an area expansion: indeed, for large $R$ we have
%\begin{eqnarray*} \textrm{Area}(\{p: R_1\le r(p)\le
%R\})&=& 2\pi \left(\ln R - \frac\eta{2R^2} + O(R^{-2-\epsilon})\right) \\
%&=& 2\pi \ln R - \frac{m-2\pi}{2R^2} + O(R^{-2-\epsilon})\;.
%\end{eqnarray*}
This suggests strongly that some methods specific to dimension
$2$, perhaps in the spirit of the Shiohama theorem ({\em cf.}\/
also~\cite{LiTam}), could provide simpler
proofs of geometric invariance and positivity when $n=2$; we have
not investigated this issue any further.
\end{remk}

\section{Mass positivity for metrics asymptotic to the standard
hyperbolic metric}
\label{S4}

 In this section we consider  metrics asymptotic to the
standard hyperbolic metric $b$ of constant negative curvature
$-1$; by this we mean that the Riemannian manifold $(N,\zh)$ is
the unit round sphere $(\cerc^{n-1},\can)$. We wish to show that
the usual positivity theorem holds under the weak asymptotic
hypotheses considered in the previous sections.

\begin{theo}\label{theo:positive} Let $(M,g)$ be a
complete boundaryless spin manifold with a $C^2$ metric, and
with scalar curvature satisfying
$$R_g \geq  -n(n-1)\;,$$ and suppose that the asymptotic
conditions \eq{m3} and \eq{m5} hold with
$(N,\zh)=(\cerc^{n-1},\can)$. Then the covector $p_{(\mu)}$
defined by \Eq{pmu} is \emph{timelike future directed or zero}
\upn{(}in particular $p_{(0)}\ge 0$\upn{)}. Moreover, it vanishes
if and only if $(M,g)$ is isometrically diffeomorphic to the
hyperbolic space.
\end{theo}

\begin{remk}
1. The $C^2$ differentiability of the metric can be replaced by a
weighted $W^{2,p}$ Sobolev condition.

\medskip

2. As already pointed out, we say that a linear functional $p$ on
$\cNb$ is \emph{causal} (resp. \emph{timelike}) and
\emph{future-directed} if it can be written as $(p_{(0)},\ldots,
p_{(n)})$ in any orthonormal and future-oriented basis
$(V_{(0)},\ldots,V_{(n)})$ with
\begin{equation}\label{remindercausal}
(p_{(0)})^2 - \sum_{i=1}^n (p_{(i)})^2 \geq 0  \
 \textrm{ and } \
p_{(0)} \geq 0 \ \textrm{ \upn{(}resp. } > 0\upn{)}.
\end{equation}
This is the obvious equivalent of the corresponding definition for
vectors; note, however, that with our signature $(-,+,\ldots,+)$
future directed vectors are  {\em not} mapped to future directed
covectors by the isomorphism of $TM$ with $T^*M$ associated with
the metric.

 We emphasise that the Lorentz vector character of
$p_{(\mu)}$ is {\em not} related to the tangent space of some
point of $M$, or of some ``abstract asymptotic point" (``the
tangent space at $i^o$" --- this last interpretation can be given
to energy-momentum in the asymptotically Euclidean context), but
arises from the fact that the adjoint action of the isometry group
of the (standard) hyperbolic space, on the subspace of its Lie
algebra singled out by Equations~\eq{eq:1}-\eq{eq:2}, is that of
the defining representation of the Lorentz group on $\R^{n+1}$.

{\flushleft 3.} Condition \eq{m3b} is actually not necessary for
positivity, in the following sense: under the remaining conditions of
Theorem~\ref{theo:positive}, the argument of the proof of
Proposition~\ref{P1} shows that $p_{(0)}=\infty$ whenever \eq{m3b}
does not hold.

{\flushleft 4.} Such a theorem cannot be obtained in a more general
setting. For
instance, in the asymptotically Euclidean context, it is
well-known that positivity statements may fail if the metric is
asymptotic to some $\Z_2$-quotient of the Euclidean space
\cite{lebrun-cex}. In the asymptotically hyperbolic setting,
Horowitz and Myers~\cite{HorowitzMyers} have constructed an
infinite family of metrics with ends asymptotic to a cuspidal
hyperbolic metric (the topology of the end is $\R\times T^2$), and
with  masses as negative as desired. If the topology of the end is
the product of a half-line with a negatively curved Riemann
surface, the mass may also be negative when minimal interior
boundaries are allowed, and it is expected that the infimum of the
possible masses is achieved only for the Kottler black-hole
metrics \cite{HorowitzMyers,Kottler,Vanzo,BLP,ChruscielSimon}.

{\flushleft 5.} The $C^2$ differentiability of the metric can be replaced
by weighted Sobolev-type conditions; this is, however, of no
concern to us here.
\end{remk}

As a corollary of Theorem~\ref{theo:positive}, together with
\cite[Theorem~I.3]{ChruscielSimon} and the remarks at the end of
Section~V of \cite{ChruscielSimon} one has (see
\cite[Corollary~I.4]{ChruscielSimon}; compare~\cite{BGH}):

\begin{theo}\label{Ctheo} Let $V$ be a strictly positive function on a three
dimensional manifold $M$ such that the metric
$$\gamma:=-V^2dt^2+g\;,$$
is a static solution of the vacuum Einstein equations with
strictly negative cosmological constant on the  space-time
${\mycal M}:=\R\times M$. If
\begin{enumerate}
\item $(M,g)$ is $C^3$ compactifiable  in the sense of
Section~\ref{Sccm} below, and if
\item the conformal boundary  at infinity of $M$
is $S^2$, with $V^{-2}g$ extending by continuity to the unit round
metric on $S^2$,
\end{enumerate}
then $(M,g)$ is the hyperbolic space, so that $({\mycal
M},\gamma)$ is the anti-de Sitter space-time.
\end{theo}

{\flushleft\it Preliminaries to the proof}. The proof of
Theorem~\ref{theo:positive} will follow the
Gibbons-Hawking-Horowitz-Perry variation \cite{GHHP} of the
classical Witten argument for the positivity of mass
\cite{witten:mass} (\emph{cf.}\/~also~\cite{min-oo:scalar,AndDahl}
and the remarks done in \cite{Delay}), and relies on the existence
on the hyperbolic space of a wealth of distinguished spinor
fields, called {\sl imaginary Killing spinors}. These are
solutions $\varphi$ of the differential equation
\begin{equation}
\widehat{\mzD}^b_X\varphi = \mzD^b_X\varphi + \frac{i}{2}
c_b(X)%\cdot
\varphi = 0\;,
\end{equation}
where we denote by $c_b(X)\varphi$  the Clifford action of a
vector $X$ on a spinor $\varphi$ with respect to the metric $b$.
On hyperbolic space there is a set of maximal dimension of
imaginary Killing spinors, which trivialise the spinor bundle.
They can be described explicitly in the following manner
\cite{baum-killing}: one may choose the standard basis $\{\p_i\}$
of the flat space as reference frame, thus inducing an isomorphism
between the spinor frame bundle of $(\B,e)$ and $\B\times
\mathit{Spin}(n)$. This can be transferred to the hyperbolic space
through the usual conformal covariance (of ``weight zero'') of
spinor bundles \cite{pg-italien}. In this trivialisation, the
Killing spinors of $b$ are then the spinor fields $\varphi_u$
given by
\begin{equation}\label{killing}
\varphi_u (x) = \omega^{-\frac{1}{2}}\left( 1 - i c_{\delta}(x) \right) u
\end{equation} where $u$ is any non-zero
constant spinor on the flat ball $\B$, and $\omega$ is the conformal
factor of the hyperbolic metric defined in \eq{omf}.

Following the terminology due to H.~Baum, Th.~Friedrich and
I.~Kath~\cite{frkath,baum-dim5,baum-killing}, the spinor
$\varphi_u$ is said to be of type I (resp. of type II) if
\begin{equation}\label{type}
\|u\|_\delta^4 + \sum_{i=1}^{n} \la c_\delta(\p_i) u , u
\ra_\delta^2 \ \textrm{ is zero (resp. is positive).}
\end{equation}
Type I spinors are actually sufficient for our purposes, we shall
describe these ones only. For any imaginary Killing spinor, the
function $$V_u = \la\varphi_u,\varphi_u\ra_b$$ is always an
element of $\cNb$. If $\varphi_u$ is moreover of type I, then
there is a set of $n$ constants $(a_i)\in \cerc^{n-1}\subset \R^n$
and a constant $\lambda>0$, such that $V_u = \lambda (V_{(0)} -
\sum_i a_i V_{(i)})$: Indeed, an explicit computation from
Equation (\ref{killing}) above shows that, in the ball model,
\begin{equation}\label{eq:V_utype1}
V_u (x) \ = \ ||u||^2_{\delta} \,\frac{1+|x|^2}{1-|x|^2} \ + \
i\sum_{j=1}^{n} \langle c_{\delta}(\partial_j)u,u\rangle_{\delta}
\, \frac{2\, x^j}{1-|x|^2}  \ , \end{equation} which is clearly
future directed, and the type I condition (\ref{type}) yields that
$V_u$ is isotropic in $\cNb$. This shows in particular that
Killing spinors of type I always exist on the hyperbolic space in
any dimension. Further, \emph{e.g.} as a result of covariance
under isometries, any future directed isotropic combination
$V_{(0)} - \sum_i a_i V_{(i)}$, $(a_i)\in \cerc^{n-1}$, can be
obtained as a $V_u$ for some type I Killing spinor ({\it i.e.} for
some constant spinor $u$ on the flat ball) in any dimension. For
later use we also note that
$$dV_u(X) = i\la c_b(X)\varphi_u,\varphi_u\ra_b\;.$$

{\flushleft\it Proof of Theorem \ref{theo:positive}}. Let $A$ be
the symmetric endomorphism defined over $\Mext$ by
$g(A\cdot,A\cdot) = b(\cdot,\cdot)$, which we will take to be of
the form $$A = I - \frac{1}{2}\, e + \{\textrm{quadratic and
higher order in } e\}$$ if $e$ is small enough; by this we mean
that $A^i{}_j=\delta^i_j - \frac 12 b^{ik}e_{kj} + $ a second
order Taylor expansion error term.  One may use $A$ as an
isomorphism between the orthonormal frame bundles of $b$ and $g$
and any lift of it as an isomorphism between their spinor frame
bundles. This enables, as in \cite{AndDahl}
(compare~\cite{Bourguignon92}), to transfer the spin connection
$\mzD^b$ of $(\Mext,b)$ on the spinor bundle of $(\Mext,g)$; for
notational convenience, the new connection will be denoted by
$\mzD^{\hat b}$. Note that this has the effect that the Clifford
action $c_b(X)$ of a vector $X$ is transformed into the Clifford
action $c_g(AX)$ of $AX$. As a consequence, the transferred
spinors, still denoted by $\varphi_u$, are now solutions of
\begin{equation}\label{eq:new-kill}
\widehat{\mzD}^{\hat b}_X\varphi_u = \mzD^{\hat b}_X\varphi_u +
\frac{i}{2} c_g(AX) \varphi_u = 0.
\end{equation}
We now denote by $\mzD$ the spinor connection associated
to the Levi-Civita connection of the metric $g$ and define the
modified connection on spinors
$$\widehat{\mzD}_X = \mzD_X + \frac{i}{2}
c_g(X).$$
 For any $\varphi_u$ we set
 $$\Phi_u =\chi\varphi_u + \psi_u\;,$$ where $\chi$ is a cut-off
function that vanishes outside of $\Mext$ and is equal to $1$ for
$r$ large enough. Suppose, first, that $\psi_u$ is compactly
supported, hence vanishes for $r\ge R$ for some $R$ on $\Mext$. We
apply the standard Schr\"odinger-Lichnerowicz
\cite{Lichnerowicz63,Schrodinger32} formula relating the rough
Laplacian of the modified connection $\widehat{\mzD}$ to the Dirac
Laplacian $\widehat{\dirac}^*\widehat{\dirac}$ \cite{AndDahl},
where
\begin{equation}
\widehat{\dirac}\Phi_u = \dirac \Phi_u - \frac{ni}{2}\Phi_u\;,
\end{equation}
with $\dirac$ being the usual Dirac operator associated with the
metric $g$. Letting $S_R=\{r=R\}\subset \Mext$ one obtains
\begin{eqnarray}\label{weitzenbock0}
&& \\\nn \int_{M\setminus\{r\ge R\}} \|\widehat{\mzD}\Phi_u \|^2_g
+ \frac{1}{4}\left(R_g + n(n-1)\right) \|\Phi_u\|_g^2 -
\|\widehat{\dirac}\Phi_u \|^2_g &=& \int_{S_R} B_{A\nu}(\Phi_u)
\\\nn&=& \int_{S_R} B_{A\nu}(\varphi_u)\;,
\end{eqnarray}
where $\nu$ is the outer $b$-unit normal to $S_r$, so that $A\nu$
is its outer $g$-unit normal, and $B_{A\nu}(\varphi_u)$ is the
boundary integrand, explicitly defined by
\begin{equation} B_{Y}(\rho) = \la \widehat{\mzD}_{Y} \rho
+ c_g(Y) \widehat{\dirac}\rho , \rho \ra_g \end{equation} for any
spinor $\rho$ and vector $Y$.

Assume that $(M,g)$ is \emph{not} the hyperbolic space, otherwise
there is nothing to prove. Let $H$ be the usual
Hilbertian completion of the space of compactly supported smooth
spinors $\psi$ on $M$ with respect to the norm defined as
\begin{equation}\label{norm}
\|\psi\|^2_H:=\int_{M}\left( \|\widehat{\mzD}\psi\|^2_g +
\frac{1}{4}\left(R_g + n(n-1)\right) \|\psi\|_g^2 \right)d\mu_g\;.
\end{equation}
We wish to show that for any $\Phi_u=\chi\varphi_u+\psi_u$, with
$\psi_u\in H$, we will have
\begin{equation}\label{weitzenbock1}
\int_{M} \|\widehat{\mzD}\Phi_u \|^2_g + \frac{1}{4}\left(R_g +
n(n-1)\right) \|\Phi_u\|_g^2 - \|\widehat{\dirac}\Phi_u \|^2_g =
\lim_{R\to\infty}\int_{S_R} B_{A\nu}(\varphi_u)\;.
\end{equation}
We start by showing that $H$ can be identified with a space $\mcH$
of $H^1_\loc$ spinor fields on $M$, with the norm $\|\cdot\|_H$
still given by \eq{norm} (after identification) for all $\psi\in
H$. First, it is not too difficult to show~\cite{BartnikChrusciel}
that in dimension larger than or equal to three there exists a
strictly positive $L^\infty_\loc$ function $w$ on $M$ such that
for all $H^1_\loc$ spinor fields $\psi$ with compact support we
have
\bel{wpimc} \int_M \|\psi\|^2_g \,w \,d\mu_g\le \int_M
\|\widehat{\mzD} \psi\|^2_g d\mu_g\;.\ee The function $w$ can be
chosen to be constant in the asymptotically hyperbolic end. In
dimension two one can also prove \eq{wpimc} if one assumes further
that there are no imaginary Killing spinors. This is sufficient
for our purposes because, if there exists a Killing spinor then,
by \cite{baum-killing}, we are in hyperbolic space, where there is
nothing to prove. So one might as well suppose that there are no
such spinors.

 Let
$\mcH$ be the space of measurable spinor fields on $M$ such that
\bel{wpimc2}\|\psi\|_{\mcH}^2:= \int_M\|\psi\|^2_g \left(w+ \frac 14
\left(R_g+n(n-1)\right)\right) d\mu_g + \int_M \|
\widehat{\mzD}\psi\|^2_g d\mu_g< \infty\;\ee where
$\widehat{\mzD}\psi$ is understood in the distributional sense.
Define $\zmcH\subset \mcH$ as the completion of $C^\infty_c$, in
$\mcH$, with respect to the $\|\cdot\|_\mcH$ norm. It is then easy
to verify the following:
\begin{Proposition} \label{Pwpimc0} The inequality
\eq{wpimc} remains true for all $\psi\in \zmcH$.
\end{Proposition}

\proof Both sides of \eq{wpimc} are continuous on
$(\mcH,\|\cdot\|_\mcH)$. \qed

\begin{Proposition} \label{Pwpimc} If $(M,g)$ is complete then $\mcH=\zmcH$.
\end{Proposition}

\proof If $\phi\in \mcH$ then the sequence $\chi_i\phi$ converges
to $\phi$ in $(\mcH,\|\cdot\|_\mcH)$, where
$\chi_i(p)=\chi(d_{p_0}(p)/i)$, where $d_{p_0}$ is the distance to
some chosen point $p_0\in M$, while $\chi:\R\to [0,1]$ is a smooth
function such that $\chi|_{[0,1]}=1$, $\chi|_{[2,\infty)}=0$.
Smoothing $\chi_i\phi$ using the usual convolution operator yields
the result. \qed

\begin{Proposition}  If $(M,g)$ is complete then there is a natural
continuous bijection between $(\mcH,\|\cdot\|_\mcH)$ and
$(H,\|\cdot\|_H)$ which is the identity on $C^1_c$; in particular,
elements of $H$ can be identified with spinor fields on $M$ which
are in $\mcH$.
\end{Proposition}

\proof By Proposition~\ref{Pwpimc} both spaces are Hilbert spaces
containing $C^1_c$ as a dense subspace, with the norms being
equivalent when restricted to $C^1_c$ by
Proposition~\ref{Pwpimc0}.\qed

Let $F(\psi)$ denote the left-hand-side of \Eq{weitzenbock1} with
$\Phi_u=\chi\varphi_u+\psi$ there, let $\psi_i\in C^1_c$ converge
to $\psi$ in $H$, we have
\begin{eqnarray*}
F(\psi)-F(\psi_i) & = & \|\psi\|^2_H-\|\psi_i\|^2_H \\
& & +2 \int_{M}
\langle\widehat{\mzD}(\chi\varphi_u),\widehat{\mzD}(\psi-\psi_i)\rangle \\
&& -  2 \int_{M}
\langle\widehat{\dirac}(\chi\varphi_u),\widehat{\dirac}(\psi-\psi_i)\rangle \\
&&  + \frac{1}{2}\int_{M}\left(R_g + n(n-1)\right)
\langle\chi\varphi_u,\psi-\psi_i\rangle\;.
\end{eqnarray*}
 It should be clear from the fact that
$\widehat{\mzD}(\chi\varphi_u)\in L^2(M)$ that all the terms above
converge to zero as $i$ tends to infinity, except perhaps for the
last one (recall that we are only assuming that $0\le (R_g+n(n-1))
|V| \in L^1(\Mext)$); the convergence of that last term can be
justified as follows:
\begin{eqnarray*}\left|\int_{M}\left(R_g + n(n-1)\right)
\langle\chi\varphi_u,\psi-\psi_i\rangle\right| &\le&
\left(\int_{M}\left(R_g + n(n-1)\right)
\|\chi\varphi_u\|^2_g\right)^{1/2}
\\
&& \times  \left(\int_{M}\left(R_g + n(n-1)\right)
\|\psi-\psi_i\|^2_g\right)^{1/2}
\\
&\le& \|\chi\varphi_u\|_H\|\psi-\psi_i\|_H\;.
\end{eqnarray*}
(Here we have applied the Cauchy-Schwarz inequality associated
with the positive quadratic form occurring in the left-hand-side
above.) Now, $F(\psi_i)=F(0)$, and we have shown that
\Eq{weitzenbock1} holds for all $\psi_u\in H$, as claimed.

 To obtain positivity of the left-hand-side of
\Eq{weitzenbock1} we seek a $\Phi_u$ such that
\begin{equation}\label{Witteneq}
\widehat{\dirac}\Phi_u = 0 \quad \Longleftrightarrow \quad
\widehat{\dirac}\psi_u=- \widehat{\dirac}(\chi\varphi_u)\;.
\end{equation}
We now use the fact that $\varphi_u$ solves (\ref{eq:new-kill}),
hence \be \widehat{\mzD}_X\varphi_u = (\mzD_X - \mzD^{\hat
b}_X)\varphi_u - \frac{i}{2}c_g(AX -X)\varphi_u\;.
\label{covphiu}\ee Now, in any $g$-orthonormal frame
$\{f_{\alpha}\}_{\alpha = 1,...,n}$, if $\omega$ denotes the
connection $1$-forms of either $\mzD$ or $\mzD^{\hat b}$ (with the
obvious notations), one has:
\be\label{covphiu1} \mzD_X - \mzD^{\hat b}_X = \frac{1}{4}
\sum_{\alpha,\beta = 1}^n ( \omega_{\alpha\beta}(X) - \omega^{\hat
b}_{\alpha\beta}(X) ) c_g(f_{\alpha})c_g(f_{\beta}) \;. \ee The
calculations of \cite[Section~2.2]{AndDahl} lead then to
\begin{equation}\label{ltwoineq1} |\widehat{\dirac}\varphi_u |\le
C\left(|\zD A|_b + |A-\textrm{id}|_b\right)\; |\varphi_u| \quad
\Longrightarrow \quad \widehat{\dirac}\left(\chi\varphi_u
\right)\in L^2(M, d\mu_g)\;.
\end{equation}
\Eq{ltwoineq1} and arguments known in principle
\cite{AndDahl,ParkerTaubes82,ChBlesHouches,Herzlich:mass} ({\em
cf\/.} also~\cite{BartnikChrusciel}) show
that there exists $\psi_u$ in $H$ such that
\Eq{Witteneq} holds. It then remains to show that the integral at
the right-hand-side of \Eq{weitzenbock1} is related to the map
$H_\Phi$. In order to do this, we complement $\nu$ into a direct
$b$-orthonormal basis $\{ \nu, e_i\}_{i=1,...,n-1}$. Seen as
sitting in $M$, $\{ A\nu, Ae_i\}_{i=1,...,n-1}$ is a direct
$g$-orthonormal basis on $S_r$. One easily finds
\begin{equation}\label{eq:termeB} B_{A\nu} (\varphi_u)\  =
\,\sum_{i=1}^{n-1} \la c_g(A\nu)
c_g(Ae_i)\widehat{\mzD}_{Ae_i}\varphi_u, \varphi_u\ra_g \;.
\end{equation}
{}From \eq{covphiu} and \eq{covphiu1}, we obtain
\begin{equation*}\begin{split}
B_{A\nu} (\varphi_u)\  =\ & \frac{1}{4}\sum_{i=1}^{n-1}
\sum_{\alpha,\beta=1}^n (\omega_{\alpha\beta}(Ae_i) - \omega^{\hat
b}_{\alpha\beta}(Ae_i) ) \la
c_g(A\nu)c_g(Ae_i)c_g(f_{\alpha})c_g(f_{\beta}) \varphi_u,
\varphi_u
\ra \\
& +\frac{i}{2} \sum_{i=1}^{n-1} \la
c_g(A\nu)c_g(Ae_i) c_g(A(Ae_i) -Ae_i)\varphi_u, \varphi_u\ra\\
= \ & - \frac{1}{2}\sum_{i,j=1}^{n-1}
(\omega_{j0}(Ae_i) - \omega^{\hat b}_{j0}(Ae_i) )
\la c_g(Ae_i)c_g(Ae_j) \varphi_u, \varphi_u \ra  \\
& + \frac{1}{4}\sum_{i,j,k=1}^{n-1}
(\omega_{jk}(Ae_i) - \omega^{\hat b}_{jk}(Ae_i) )
\la c_g(A\nu)c_g(Ae_i)c_g(Ae_j)c_g(Ae_k) \varphi_u, \varphi_u
\ra \\
& +\frac{i}{2} \sum_{i=1}^{n-1} \la c_g(A\nu)c_g(Ae_i) c_g(A(Ae_i)
-Ae_i)\varphi_u, \varphi_u\ra\;,
\end{split}\end{equation*}
where the subscript $._0$ in the last formula denotes the basis
element $A\nu$. These formulae correct Equation (34), page 20,
in~\cite{AndDahl}: in the second line of that equation the
multiplicative factor 1/4 should be changed to 1/8; this minor mistake
carries over to all the equations that follow.

One may now use \cite[Formulas (2--3)]{AndDahl} to compute the
difference between the connection $1$-forms of $\mzD$ and $\mzD^b$
(or, equivalently, the connection $1$-form of $\mzD^{\hat b}$)
with respect to the covariant derivative $\mzD^bA$. Following
again Andersson and Dahl's argument \cite{AndDahl}, and noting
that all the imaginary-valued terms have to cancel out
 because the
left-hand-side of \eq{weitzenbock0} is real, one eventually gets:
\begin{equation}\label{eq:etape1}\begin{split} B_{A\nu}& (\varphi_u)\,
 = \ \frac{1}{2}\,\sum_{i=1}^{n-1} \left( g((\mzD^b_{A\nu}A)e_i,
Ae_i)-
g((\mzD^b_{Ae_i}A)\nu, Ae_i) \right) \la\varphi_u,\varphi_u\ra_g \\
& \ \ \ + \frac{1}{4}\,\sum_{i,j,k \textrm{ distinct}} g(
(\mzD^b_{e_i}A)e_j, Ae_k)\,\la c_g(A\nu)c_g(Ae_i)c_g(Ae_j)c_g(Ae_k)
\varphi_u ,\varphi_u\ra_g \\
& \ \ \ + \frac{i}{2} \sum_{i=1}^{n-1} \,\la
c_g(A\nu)c_g(Ae_i)c_g\left(A(Ae_i) - Ae_i\right)
\varphi_u,\varphi_u\ra_g.
\end{split}\end{equation}
Using the (spin) isomorphism $A$, this can immediately be
rewritten as:
\begin{equation}\label{eq:etape2}\begin{split} B_{A\nu}(\varphi_u)\
= & \ \frac{1}{2}\,\sum_{i=1}^{n-1} \left(
b(A^{-1}(\mzD^b_{A\nu}A)e_i,
e_i)-b(A^{-1}(\mzD^b_{Ae_i}A)\nu, e_i) \right)
\la\varphi_u,\varphi_u\ra_b \\
& + \frac{1}{4}\,\sum_{i,j,k \textrm{ distinct}}
b(A^{-1}(\mzD^b_{e_i}A) e_j,e_k)
\,\la c_b(\nu)c_b(e_i)c_b(e_j)c_b(e_k) \varphi_u, \varphi_u\ra_b \\
& + \frac{i}{2}\sum_{i=1}^{n-1}\,\la c_b(\nu)c_b(e_i)c_b\left(Ae_i
- e_i\right) \varphi_u,\varphi_u\ra_b ,
\end{split}\end{equation}
where all computations take now place on the spinor bundle of the
reference hyperbolic metric. Taking into account the relationship
between the squared norm of $\varphi_u$ and $V_u$, we now recall
that our asymptotic conditions imply that, in the last formula,
any quadratic term in $A-I$ and $\mzD^b A$ (or, equivalently, in
$e=g-b$), when integrated on $S_r$, has limit value zero as $r$
goes to infinity. One may then eliminate a large number of
occurrences of the map $A$ from the above formula. We will use
below the notation $ U \simeq V$ to mean that $V$ is the only term
that contributes when integrating $U$ on larger and larger
spheres. Equivalently,
$$ U \simeq V \quad\Longrightarrow \quad \lim_{r\to\infty} \int_{S_r} U
\ = \ \lim_{r\to\infty} \int_{S_r} V . $$
Then, one obtains in our case:
\begin{equation}\label{eq:etape3}\begin{split}
B_{A\nu}(\varphi_u)\ \simeq & \
 \frac{1}{2}\,\sum_{i=1}^{n-1} \left( b((\mzD^b_{\nu}A)e_i, e_i)-
b((\mzD^b_{e_i}A)\nu, e_i) \right) \la\varphi_u,\varphi_u\ra_b \\
& + \frac{1}{4}\,\sum_{i,j,k \textrm{ distinct}}
b((\mzD^b_{e_i}A)e_j,e_k)
\,\la c_b(\nu)c_b(e_i)c_b(e_j)c_b(e_k) \varphi_u, \varphi_u\ra_b\\
& + \frac{i}{2}\sum_{i=1}^{n-1}\,\la c_b(\nu)c_b(e_i)c_b\left(Ae_i
- e_i\right) \varphi_u,\varphi_u\ra_b.
\end{split}\end{equation}
Furthermore, the last term in the previous formula is easily computed
and it remains:
\begin{equation*}\begin{split}
B_{A\nu}(\varphi_u)\ \simeq & \
 \frac{1}{2}\,\sum_{i=1}^{n-1} \left( b((\mzD^b_{\nu}A)e_i, e_i)-
b((\mzD^b_{e_i}A)\nu, e_i) \right) \la\varphi_u,\varphi_u\ra_b \\
& + \frac{1}{4}\,\sum_{i,j,k \textrm{ distinct}}
b((\mzD^b_{e_i}A)e_j,e_k)
\,\la c_b(\nu)c_b(e_i)c_b(e_j)c_b(e_k) \varphi_u, \varphi_u\ra_b\\
& - \frac{i}{2}\left(\tr(A-I)\,\la
c_b(\nu)\varphi_u,\varphi_u\ra_b - \la
c_b((A-I)\nu)\varphi_u,\varphi_u\ra_b\right).
\end{split}\end{equation*}
The second term in this last formula  contributes as zero, since
$A$ (hence $\mzD^bA$) is symmetric and $e_i\cdot e_j\cdot{}$ is
antisymmetric. We then get:
\begin{equation*} \begin{split}
\lim_{r\to\infty} \int_{S_r} B_{A\nu}(\varphi_u)\ = & \
\lim_{r\to\infty}\int_{S_r} \frac{1}{2}\,\sum_{i=1}^{n-1}\left(
b((\mzD^b_{\nu}A)e_i, e_i)-
b((\mzD^b_{e_i}A)\nu, e_i) \right) \la\varphi_u,\varphi_u\ra_b \\
& - \frac{i}{2}\left(\tr(A-I)\la c_b(\nu)\varphi_u,\varphi_u\ra_b
- \la c_b((A-I)\nu)\varphi_u,\varphi_u\ra_b\right).
\end{split}\end{equation*}
To compare with our previous formulas for mass integrals, we now
replace $A$ by $I -\frac{1}{2}\, e$ (as remainder terms in a
Taylor expansion contribute as zero in the limit) and relate the
norms $\la\varphi_u,\varphi_u\ra_b$ and $\la
c_b(\nu)\cdot\varphi_u,\varphi_u\ra_b$ to $V_u$ using
$|\varphi_u|^2_b = V_u$ and $i\la c_b(X)\varphi_u, \varphi_u\ra_b
= dV_u(X)$ to obtain in the limit:
\begin{equation*}
- \frac{1}{4}\, \lim_{r\rightarrow\infty} \int_{S_r}
   V_u\, (\, d(\tr_b\, e)(\nu) - \sum_{i=1}^{n-1} \mzD^b_{e_i}e(\nu,e_i)\, )
\ - \ (\tr_b\, e)\, dV_u(\nu) + dV_u (e(\nu))\ ,
\end{equation*}
or, equivalently,
\begin{equation}
\lim_{r\to\infty} \int_{S_r} B_{A\nu}(\varphi_u)\ = \ \frac{1}{4}\,
\lim_{r\rightarrow\infty} \int_{S_r} \ourU^i\nu_i.
\end{equation}
As the left-hand side in the Schr\"odinger-Lichnerowicz formula
(\ref{weitzenbock1}) is nonnegative, this implies that the mass
(seen as a linear functional on $\cNbz$) is non-negative on any
future-directed null vector in $\cNbz$. Standard considerations in
Lorentzian geometry
yield that this linear functional is causal and future directed
(notice that we use here only the existence of imaginary Killing
spinors of type I, which is valid in any dimension). It remains to
show that it is either timelike, or vanishes.

Suppose, then, that the mass $H_\Phi$ (still seen as a linear
functional on $\cNbz$) is isotropic, then there exists a non-zero
element $W$ of the light cone in $\cNbz$ such that it evaluates
against $W$ as zero. Up to rescaling, $W$ can be written as
$V_{(0)} - \sum_i a_i V_{(i)}$ and we already noticed there exists
a constant spinor $u$ on $(\B,e)$ such that $V_u=W$. The
Lichnerowicz-Schr\"odinger formula above has then a vanishing
contribution at infinity. This implies the associated spinor
$\Phi_u$ is a Killing spinor, and H.~Baum's work shows that
$(M,g)$ is isometric to the hyperbolic space
\cite{baum-killing,AndDahl}. The mass functional $H_\Phi$ is then
zero and this ends the proof of Theorem \ref{theo:positive}.\qed

\medskip

{\flushleft\it Manifolds with boundary}. When $(M,g)$ has a
compact boundary one expects that the correct statement is the
Penrose inequality~\cite{HI1,HI2,Bray:preparation2}, which seems
to lie outside of the scope of the Witten-type argument given
above. Recall, however, that this last argument does lead to a
positivity statement
\cite{GHHP,Herzlich:mass,mh-inegalite-penrose} when compact
boundaries occur:

\begin{theo}\label{Tmb} Let $(M,g)$ be a
complete spin manifold with a $C^2$ metric, with a compact
non-empty boundary of mean curvature $$\Theta \leq n-1 ,$$ and with scalar
curvature satisfying
$$R_g \geq  -n(n-1)\;.$$ If the asymptotic
conditions \eq{m3} and \eq{m5} hold with
$(N,\zh)=(\cerc^{n-1},\can)$, then the covector $p_{(\mu)}$
defined by \Eq{pmu} is \emph{timelike future directed}.
\end{theo}

\begin{proof} When $\partial M$ is non-empty, a supplementary boundary integral
over $\partial M$, given by
\bel{supbi} \int_{\partial M}B_{A\nu}(\Phi_u) = \int_{\partial M}
\langle\dirac_{\partial M} \Phi_u + \frac 12
\left(\Theta-(n-1)ic_g(n)\right)\Phi_u,\Phi_u\rangle\;, \ee
appears in \Eqsone{weitzenbock0} and \eq{weitzenbock1}, where
$\Theta$ is the inwards oriented mean extrinsic curvature of
$\partial M$, while $\dirac_{\partial M}$ is a boundary Dirac
operator. It is defined as $$\dirac_{\partial M} =
c_g(n)\sum_{i=1}^{n-1} c_g(e_i){\mycal D} _{e_i}\;,$$ where
$\mycal D$ is the spin connection intrinsic to $\partial M$,
explicitly defined on spinors fields on $M$ restricted to
$\partial M$ as
\begin{equation}\label{lindcon}
\mycal D_X = D_X - \frac{1}{2}c_g(n)c_g(B(X)) \;,
\end{equation}
where $B$ is the shape operator of the boundary.

Following~\cite{GHHP} (compare~\cite{Herzlich:mass}) we impose
%Atiyah-Patodi-Singer~\cite{APS75,Herzlich:mass,BartnikChrusciel}
the following boundary condition on the spinors $\Phi_u$:
\bel{bcpu} \Phi_u=\varepsilon \Phi_u \;,\ee
where $ \varepsilon$ is a hermitian involution on spinors given by
\begin{equation}\label{choice}
\varepsilon = i \, c_g(n)
\end{equation}
as in \cite{HMZ-asian}. It is proved in this paper that this leads
to a self-adjoint elliptic problem for the Dirac operator which
can be solved. Positivity of the mass is obtained through the same
argument as before, the boundary contribution \eq{supbi} having
the correct sign since $\varepsilon\dirac_{\partial
M}=-\dirac_{\partial M}\varepsilon$ so that
\begin{eqnarray}\nonumber \langle \Phi_u,\dirac_{\partial M}\Phi_u\rangle &=& \langle
\Phi_u,\dirac_{\partial M}\varepsilon\Phi_u\rangle \\\nonumber  &
=& -\langle \Phi_u,\varepsilon \dirac_{\partial M}\Phi_u\rangle
\\\nonumber
&=&-\langle \varepsilon\Phi_u, \dirac_{\partial M}\Phi_u\rangle \\
&=&-\langle \Phi_u, \dirac_{\partial M}\Phi_u\rangle
\;.\label{bvas}\end{eqnarray} As a result, $\langle
\Phi_u,\dirac_{\partial M}\Phi_u\rangle$ vanishes and it remains
\bel{supbi1} \int_{\partial M}B_{A\nu}(\Phi_u) = \int_{\partial M}
\langle \frac 12
\left(\Theta-(n-1)\right)\Phi_u,\Phi_u\rangle\; \ee
for the boundary contribution. This proves as above that the covector
$p_{(\mu)}$ defined by \Eq{pmu} is timelike future directed, or lightlike future
directed, or vanishing. Let us show that those last two
possibilities cannot occur: clearly, $p_{(\mu)}$ can be lightlike
or vanish if and only if $M$ carries an imaginary Killing spinor
$\Phi_u$ satisfying the boundary condition \eq{bcpu} at $\partial
M$. Further, Equation~(\ref{supbi1}) implies that
$\Theta$ is identically equal to $(n-1)$ ---otherwise, the imaginary Killing
spinor field $\Phi_u$ would be zero on an open set on the
boundary, a situation which is forbidden by the uniqueness
property of solutions of ordinary differential equations.
Moreover it is a classical fact \cite{baum-killing} that existence of an
imaginary Killing spinor implies that $(M,g)$ is Einstein, of scalar
curvature $-n(n-1)$.

Hence we have the following situation : a non compact Einstein
manifold, looking like the hyperbolic space at infinity and with a
compact inner boundary of constant mean curvature $n-1$. Choose
any very large sphere-like compact submanifold $S$ in the
asymptotically hyperbolic end of $M$ and consider the
part of $M$ located \emph{inside} $S$. It is a compact Einstein
$n$-manifold with boundary having two components, one of constant
mean curvature $\Theta_{\partial M} = n-1$
and the other one having (not necessarily constant)
mean curvature $\Theta_S$ close to that of a sphere in the hyperbolic space by
\eq{m5}, hence which can
be taken so that $\Theta_S > n-1$ at each point of $S$ (here, both mean curvatures
are computed with respect to the normal unit vector pointing towards infinity).

Let now $p$ be in $\partial M$ the closest point to $S$ and
$\gamma$ be a minimising geodesic from $S$ to $p$, starting from a
point $q$ in $S$. Let $\ell$ be the distance from $\partial M$ to
$S$, {\it i.e.} the distance from $p$ to $q$. We now consider the
family $\partial M_{\delta}$ of submanifolds obtained by pushing
the boundary $\partial M$ a distance $\delta$ through its normal
exponential map towards $S$, and the analogously defined
submanifolds $S_{\eta}$ obtained from $S$ by pushing it a distance
$\eta$ towards $\partial M$.

For $\delta>0$ small enough, the submanifold $\partial M_{\delta}$
is still smooth. Moreover, the contact point $r$ of $\partial
M_{\delta}$ and $\gamma$ is necessarily the closest point of
$\partial M_{\delta}$ to $S$ (and is at distance $\ell -\delta$).
As $\gamma$ is minimising, the distance function to $S$ is smooth
in an open neighborhood of $\gamma\setminus\{p\}$, hence the
submanifold $S_{\ell-\delta}$ is smooth around $r$, contained in
the (closure of the) unbounded part of $M$ delimited by $\partial
M_{\delta}$ and is necessarily tangent at $r$ to $\partial
M_{\delta}$.

It remains to show that this leads to a contradiction. This
follows from classical comparison geometry: the usual Riccati
equation for the normalised mean curvature $H=\frac{\Theta}{n-1}$
reads \cite{cheeger-LNM}:
\[ H' \leq - H^2 - \frac{\Ric(\gamma',\gamma')}{n-1} \]
where $H$ stands either for the \emph{outwards} normalised mean
curvature of the family $\{\partial M_{\delta}\}$ or for the
\emph{outwards} normalised mean curvature of the family
$\{S_{\eta}\}$ and $'$ denotes differentiation with respect to
either $\delta$ or $-\eta$. As $(M,g)$ is Einstein, this
translates in our context as:
\[ H' \leq 1 - H^2 \]
and one gets by standard arguments that
\[ H_{\partial M_{\delta}} \leq 1 \ \textrm{ and } \ H_{S_{\eta}} > 1 \
\textrm{ for any } \ \eta >0, \ \delta > 0.\] At the point $r$,
this contradicts the comparison principle for the mean curvature
equation, which ends the proof.
\end{proof}

\begin{remk} This result is of special interest in general relativity, where the
condition on the mean curvature ($\Theta\leq n-1$) has the
following interpretation: let $\alpha\in\R$ be a constant
satisfying $|\alpha|\le 1$, then our Riemannian manifold can be
thought of as arising from a spacelike slice $i(M)$ in a vacuum
space-time with cosmological constant
$\Lambda:=-(1-\alpha^2)n(n-1)\le 0$, such that $i(M)$ has
extrinsic curvature $K_{ij}=\alpha g_{ij}$. The condition $R_g\ge
-n(n-1)$ is then equivalent to requiring positive energy density
on $i(M)$, while the condition $\Theta\leq -\alpha(n-1)$ is
equivalent to the statement that $\partial M$ is an
outer-future-trapped, or marginally outer-future-trapped, compact
hypersurface in $i(M)$. Under suitable global conditions existence
of such surfaces implies existence of a black hole region in the
associated space-time. Similarly the condition $\Theta\leq
\alpha(n-1)$ is associated with outer-past-trapped surfaces, and
leads to existence of white hole regions. A significant
consequence of the above result is then that the trapped surface
situation is \emph{far away} from the case of vanishing mass. This
gives  mathematical support (disjoint from all known physical
reasons) to the idea that some statement analogous to the Penrose
inequality \cite{Bray:preparation2,HI2} should hold in the
asymptotically hyperbolic case as well.
\end{remk}

\begin{remk} In the special (and interesting for physics) case $\alpha=0$
(\emph{i.e.} $\Theta \le 0$), there is a shorter way to prove that
mass cannot vanish: in the proof of Theorem \ref{Tmb}, one may
take for $\varepsilon$ any self-adjoint involution satisfying
$c_g(n)\varepsilon= -\varepsilon c_g(n)$ and
$\varepsilon\dirac_{\partial M}=-\dirac_{\partial M}\varepsilon$
(such an $\varepsilon$ will certainly exist if our Riemannian
manifold is isometrically embedded as a Riemannian slice in a
Lorentzian $(n+1)$-dimensional manifold, or more generally, if the
spinor bundle carries a representation of the Clifford algebra of
the Lorentzian metric $\gamma = -e_0\otimes e_0+g$; in any of
those cases one sets $\varepsilon = c_\gamma(e_0)c_\gamma(n)$
---note however that, in the rest of the proof as well as in the
other parts of the paper, our discussion will stay purely
Riemannian, as opposed to
\cite{GHHP,Herzlich:mass,min-oo:scalar}). A calculation identical
to \eq{bvas} shows that the boundary integral will have the right
sign and the proof goes through without modifications, implying
that the covector $p_{(\mu)}$ is timelike future directed, or
lightlike future directed, or vanishing. Assuming one of the last
two conclusions, the equality case in the
Lichnerowicz-Weitzenb\"ock formula yields again existence of an
imaginary Killing spinor. When restricted to the boundary, this
spinor field would then be an eigenspinor of the formally
 self-adjoint boundary Dirac
operator $\dirac_{\partial M}$ for a \emph{purely imaginary}
eigenvalue, which is certainly impossible on a compact manifold.
\end{remk}

\section{The mass of conformally compactifiable asymptotically
hyperbolic ends}  \label{Sccm}
The metric $g$ of a Riemannian manifold $(M,g)$ will be said to be
{\em $C^k$ compactifiable} if there exists a compact Riemannian
manifold with boundary $(\bM\approx M\cup
\piM\cup
\partial M ,\tg)$, where
$\partial\bM=\partial M\cup\piM$ is the metric boundary of
$(\bM,\tg)$, with  $\partial M$ --- the metric boundary of
$(M,g)$, together with a diffeomorphism
$$\psi: \Int \bM\to M$$ such that \be \label{ccond} \psi^*g=
\Omega^{-2} \tg\;,\ee where $\Omega$ is a defining function for
$\piM$ ({\em i.e.,\/} $\Omega\ge 0$, $\{\Omega=0\}=\piM$, and
$d\Omega$ is nowhere vanishing on $\piM$), with $\tg$ --- a metric
which is $C^k$ up--to--boundary on $\bM$. The triple
$(\bM,\tg,\Omega)$ will then be called a \emph{$C^k$ conformal
completion} of $(M,g)$. Clearly the definition allows $M$ to have
a usual compact boundary. $(M,g)$ will be said to have a
\emph{conformally compactifiable end $\Mext$} if $M$ contains an
open submanifold $\Mext$ (of the same dimension that $M$) such
that $(\Mext,g|_{\Mext})$ is conformally compactifiable, with a
connected conformal boundary $\partial_\infty
\Mext$.

In the remainder of this section we shall assume for simplicity that
the conformally
rescaled metric $\tg$ is smooth up to boundary; it should be clear
how the conditions here can be adapted to a weighted H\"older or
Sobolev setting to allow lower differentiability compactifications
consistent with the requirements of Theorem~\ref{Tinv}.

It is easily seen, using the transformation properties of the
Riemann tensor under conformal transformations ({\em cf.,
e.g.},\/~\cite{Kuehnel})
 that for smoothly %$C^2$
compactifiable metrics all the sectional curvatures $\kappa$ of
$g$ satisfy \be\label{cm1h}
(\kappa+|d\Omega|^2_{\tg^{\#}})(p)\to_{p\to\piM}0\;,
 \ee
where $|\cdot|_k$ denotes the norm of a tensor with respect to a
metric $k$; recall that $g^{\#}$ is the metric on $T^*M$
associated to $g$.

Now, \Eq{ccond} determines only the conformal class $[\tg]$ of
$\tg$. Without loss of generality we can restrict the
representative $\tg$ of $[\tg]$ so that the metric $h_0$ induced
by $\tg$ on $\piM$ has constant scalar curvature normalised as in
\eq{cm1.1}, and  this restriction will be made in what follows.

 A compactifiable
metric will be called {\em asymptotically hyperbolic\/} in an end
$\Mext$ if \be\label{cm2h}  \forall\
{p\in\piMx}\qquad|d\Omega|^2_{\tg^{\#}}(p)=1\;.
 \ee
 In what follows we restrict our considerations to a single end $\Mext$,
 replacing $M$ by $\Mext$ we will assume that $M=\Mext$, so that
 $\piM=\piMx$.
 Whenever \eq{cm2h} holds on $\piM$, a preferred representative of
$[\tg]$ in a neighborhood of $\piM$ can be chosen by requiring
that
 \be
 \label{cm3}
|d\Omega|^2_{\tg^{\#}}\equiv 1\;.
 \ee
Using $x:=\Omega$ as the first coordinate, a coordinate system can
be constructed (in some perhaps smaller neighborhood of $\piM$) in
which $g$ takes the form $$ g= x^{-2}\left(dx^2+h_x\right)\;,$$
where $h_x$ is an $x$-dependent family of metrics on $N:=\piM$. We
 define the reference metric $b$ as \be\label{c4}
b:=x^{-2}\left(dx^2+ (1-kx^2/4)^2h_0\right)%\;,
\ee(compare \Eq{cm5}). To make contact with Section~\ref{Sminv},
we assume that  $b$ is one of the metrics considered there. If $r$
is defined by \Eq{radcoord}, then the asymptotic conditions of
Proposition~\ref{P1} and of Theorem~\ref{Tinv} will hold if and
only if
\begin{eqnarray}
 \label{c1}& 0\le i \le \lfloor
n/2 \rfloor \qquad \partial_x^i \left(h_x-
(1-kx^2/4)^2h_0\right)\Big|_{x=0}=0 & \\ \nn&
\Longleftrightarrow\qquad  h_x=(1-kx^2/4)^2h_0+o(x^{\lfloor
n/2\rfloor})\;, &\end{eqnarray}
 where
$\lfloor \alpha\rfloor$ denotes the integer value of $\alpha$, and
if
\begin{eqnarray}
 \label{c2}& R_g+n(n-1)=O(x^{
n-1})\;. &\end{eqnarray} For instance, in the physically
significant case $n=3$,  \Eq{c1} is equivalent to the requirement
that the second fundamental form of $\piM$ vanishes in the
conformal gauge \eq{cm3}.

Under \Eqs{c1}{c2}, given some compactification
$(\overline{M}_1,g_1,\Omega_1)$ of an end $(\Mext,g)$ of $M$, we
use the background \eq{c4} to define its mass, whenever the
resulting background is one of those discussed in
Section~\ref{Sminv}. (As already pointed out, when
$(\piM,h_0)=(\cerc^{n-1},\can)$ this definition coincides with
that of~\cite{Wang}.) Consider a second compactification
$(\overline{ M}_2,g_2,\Omega_2)$ of $(M,g)$ satisfying the above
requirements; it is far from clear that the resulting mass will be
the same. This turns out to be the case:

\begin{theo}\label{Tmcc}  Suppose that $(M,g)$ contains a conformally
compactifiable end
$\Mext$ such that  $(\piMx, \tg|_{\piMx})$ is one of the manifolds
considered in Section~\ref{Sminv} (that is, the scalar curvature
of $\zh$ is as in \eq{cm1.1} and either $\zh$ has strictly
negative Ricci curvature, or is flat, or is the round sphere, or a
quotient thereof). Assume, moreover, that \Eqs{c1}{c2} hold.
Then the mass of $\Mext$, as defined above,
%in Section~\ref{Sccm},
is independent of the compactification of
$\Mext$ chosen for its calculation.
\end{theo}

\proof
% \noindent\emph{Proof of Theorems~\ref{Tucc} and \ref{Tmcc}:}
As already pointed out above, we can  modify the $g_a$'s and
$\Omega_a$'s so that \be\label{cf1} |d\Omega_a|_{g^{\#}_a}=1 \ee
in a neighborhood of $\piM _a$. Using the $\Omega_a$'s as the
first coordinate, in neighborhoods of respective boundaries we can
write the metrics $g_a$ as $$ g_a = d\Omega_a^2 + h^a\;, $$ where
$h^a$ is the metric induced by $g_a$ on the level sets of
$\Omega_a$. We introduce radial coordinates $r_a$ as in
\eq{radcoord},
$$r_a =
\frac{1-k\Omega_a^2/4}{\Omega_a}\;,$$ so that
$$g = \frac{dr_1^2}{r_1^2+k}+r_1^2h^1_{AB}dv_1^A dv_1^B =
\frac{dr_2^2}{r_2^2+k}+r_2^2h^2_{AB}dv_2^Adv^B_2 \;,
$$
where $k$ is defined by \Eq{cm1.1} using the boundary metric
arising out from $g_1$, and we have denoted by $(r_a,v_a^A)$,
$a=1,2$, the corresponding local coordinates near $\piM _a$. It
follows from~\cite[Theorem~3.3]{ChNagyATMP} that the map
$$(r_1,v_1^A) \to (r_2,v_2^A) $$ extends by continuity to a differentiable map
from $\bM_1$ to $\bM_2$. Equivalently, $\phi_1^{-1}\circ\phi_2$
extends by continuity to a continuous map from $\bM_1$ to $\bM_2$.
Further, point 1 of \cite[Theorem~3.3]{ChNagyATMP}  shows that the
limit
$$   \lim_{r_1\to\infty}v_2^A(r_1,v_1^A)  $$
exists, and defines a $C^\infty$ conformal diffeomorphism $\Psi$
from $(\piM _1,h^1|_{\piM _1})$ to $(\piM _2,h^2|_{\piM_2})$:
\be\label{confeq} \Psi^*h^2|_{\piM _2} = e^\psi h^1|_{\piM_1}\;.
\ee Replacing $g_1$ by $e^\psi g_1$ and $\Omega_1$ by
$\Omega_1e^{\psi/2}$ , where, by an abuse of notation, we use the
same symbol $e^\psi$ to denote the extension of $e^\psi$ from
$\piM _1$ to $\bM_1$ such that
$$|d(\Omega_1e^{\psi/2})|_{(e^{\psi}g)^{\#}}=1\;,$$ we obtain
\Eq{confeq} with $\psi=0$, hence $h^1|_{\piM _1}$ is isometric to
$ h^2|_{\piM _2}$. As a result, Theorem~\ref{Tinv} together with
the discussion
of Section~\ref{Sminv} establishes Theorem~\ref{Tmcc}. \qed

\medskip

We note that the argument just given also proves the following:

\begin{pro}
\label{pscc} Consider $(\Mext,b)$ with a metric of the form
\eq{cm1}. Then {for every
   conformal  isometry $\Psi$ of $(N,\zh)$ there exists  }
   {$R_*>0$ and a $b$-isometric map $\Phi:
[R_*,\infty)\times N \to [R,\infty)\times N$, such that}
$$\lim_{r\to\infty} \Phi(r,\cdot) = \Psi(\cdot)\;.$$
\end{pro}

\section{Uniqueness of conformal completions}
\label{S5}

It should be clear that the invariance of the mass is related to
the question of uniqueness of conformal compactifications.  There
are several issues to address here:  $g_1$ is conformal to an
appropriate pull-back of $g_2$ on the interior of $\overline{
M}_1$, but the relative conformal factor could blow up as one
approaches the boundary of $M_1$. Even if the relative conformal
factor remains uniformly bounded both from above and below, it is
not clear whether or not it extends differentiably
--- or even just continuously ---
to the boundary. Let us show that things behave as expected, so
that conformal completions are conformally diffeomorphic in the
sense of \emph{manifolds with boundary}; notations are the same as
in the previous section.

\begin{theo}\label{Tucc}
Let $(M,g)$ be a Riemannian manifold endowed with two
$C^\infty$-conformal compactifications $(\bM_1,g_1,\Omega_1)$ and
$(\bM_2,g_2,\Omega_2)$. Then
$$\phi_1^{-1}\circ\phi_2: \Int M_2\to \Int M_1 $$
extends by continuity to a $C^\infty$ conformal up-to-boundary
diffeomorphism from $(\bM_2,g_2)$ to $(\bM_1,g_1)$, in particular
$\bM_1$ and $\bM_2$ are diffeomorphic as manifolds with boundary.
\end{theo}

\begin{remk} Results about completions of finite
differentiability can be obtained in a similar way, by chasing the
order of differentiability through various steps of the arguments
below.
\end{remk}

\proof
Let $\varphi_2$ be defined on $\Int \bM_1$ by the equation
$$\varphi_2:= \frac{\Omega_2\circ
\phi_2^{-1}\circ\phi_1}{\Omega_1}>0\;;$$ \Eq{cf1} gives
\begin{eqnarray}\label{cf2}
1 &=& |d\Omega_2|_{g^{\#}_2}^2
\\
\nn &=& |\varphi_2 d\Omega_1+\Omega_1d\varphi_2|_{g^{\#}_2}^2
\\ \nn &=&
\varphi_2 ^2|d\Omega_1|_{g^{\#}_2}^2  + 2
\varphi_2\Omega_1g^{\#}_2(d\varphi_2,d\Omega_1)
+\Omega_1^2|d\varphi_2|_{g^{\#}_2}^2
\\ \nn &=&
|d\Omega_1|_{g^{\#}_1}^2  + 2
\varphi_2^{-1}\Omega_1g^{\#}_1(d\varphi_2,d\Omega_1)
+\varphi_2^{-2}\Omega_1^2|d\varphi_2|_{g^{\#}_1}^2\;,
\end{eqnarray}
hence \be\label{cf3}
 2
\Omega_1g^{\#}_1(d(\ln\varphi_2),d\Omega_1)
=-\Omega_1^2|d(\ln\varphi_2)|_{g^{\#}_1}^2\le 0\;. \ee We can
identify  a neighborhood of $\pM_1$ with $\pM_1\times[0,x_0]$
using the flow of $g_1^{\#}(d\Omega_1, \cdot)$. \Eq{cf3} shows
that $\ln \varphi_2$ is monotonously increasing along the integral
curves of the vector field $g^{\#}_1(d\Omega_1,\cdot)$ when
$\Omega_1$ decreases, so that
%the limit
%\be \psi_2:=\lim_{\Omega_1\to0}\ln\varphi_2:\pM_1\to
%\R\cup\{-\infty\} \;,\ee taken along those integral curves,
%exists. By monotonicity,
there exists a constant $C_2:=\inf_{\pM_1\times\{x_0\}}\varphi_2$
such that on $\pM_1\times[0,x_0]$ we have\be\label{cf5} \varphi_2
\ge C_2>-\infty \;.\ee
 Applying the same argument, with $g_1$ and
$g_2$
interchanged, to $$ %\frac 1
{\varphi_1}:= \frac{\Omega_1\circ
\phi_1^{-1}\circ\phi_2}{\Omega_2}=\frac{
1}{\varphi_2}\circ\phi_1^{-1}\circ\phi_1: \Int \bM_2\to \R\;,$$
shows that
%gives the existence of the limit \be \psi_1:=
%\lim_{\Omega_2\to0}\{-\ln(\varphi_2\circ\phi_2^{-1}\circ\phi_1)\}:\pM_2\to
%\R\cup\{-\infty\} \;,\ee satisfying
on $\pM_1\times[0,x_0]$ it holds
 \be\label{cf6}  \varphi_1
\ge C_1>-\infty \;.\ee Equations~\eq{cf5} and \eq{cf6} clearly
imply that the $\varphi_a$'s are uniformly bounded and uniformly
bounded away from zero.

 Set
$$\phi_{12}:=\phi_1^{-1}\circ\phi_2\;, \quad
\phi_{21}:=\phi_2^{-1}\circ\phi_1\;.$$ Let $\dga$ denote the
distance function associated with the metric $g_a$. For $p,q$ in
$\Int\bM_1$ we have
\begin{eqnarray}\label{lipcont}
\dgt(\phto(p),\phto(q))&=&\inf_ \Gamma\int_\Gamma \sqrt{g_2(\doto
\Gamma, \doto \Gamma)}
\\ &= &
\inf_ \Gamma\int_\Gamma \sqrt{\varphi_1^{-2} \,
(\phot{}^*g_1)(\doto \Gamma, \doto \Gamma)}\nn
\\ &\ge &
\inf_ \Gamma C\int_\Gamma \sqrt{ \, (\phot{}^*g_1)(\doto \Gamma,
\doto \Gamma)}\nn
\\ &=&
C\inf_{ \phot(\Gamma)} \int_{\phot(\Gamma)} \sqrt{ g_1(\doto
{\phot{}(\Gamma)}, \doto{ \phot{}(\Gamma)})} \nn
\\ \nn&=& C \dgo(p,q)\;.
\end{eqnarray} This, together with an identical calculation with
$g_1$ and $g_2$ interchanged shows that $\phot$ and $\phto$ are
uniformly Lipschitz continuous.

Clearly, $\bM_2$ is the metric completion of the manifold $M$ with
respect to the metric $(\phi_2^{-1}){}^*g_2$; similarly for
$\bM_1$. An identical calculation shows that the metrics
$(\phi_a^{-1}){}^*g_a$, $a=1,2$ define uniformly equivalent
distance functions. But completions obtained using equivalent
distances are homeomorphic; it follows that $\bM_1$ is
homeomorphic to $\bM_2$, in particular $\pM_1$ is homeomorphic to
$\pM_2$. In fact, by definition we have \bel{ideq} \phto\circ\phot
= \textrm{id}_{M_2}\;,\quad \phot\circ\phto =
\textrm{id}_{M_1}\;.\ee Since $\phot$ and $\phto$ are continuous,
they have an extension by continuity to the metric completed
spaces; we will use the same symbol to denote those extensions. It
is then easily seen that \eq{ideq} with $M_a$ replaced by $\bM_a$
holds for the extensions, so that the extensions do directly
provide the desired homeomorphism. \Eq{lipcont}, together with its
equivalent with $g_1$ interchanged with $g_2$, further show that
the extensions $\phto$ and $\phot$ are uniformly Lipschitz
continuous on $\bM_1$ and $\bM_2$. Obviously
$$\phto: \pM_1\to\pM_2\;, \quad \phot:\pM_2\to\pM_1\;,$$
with $\phto|_{\pM_1}$ and $\phot|_{\pM_2}$ being homeomorphisms
inverse to each other by the completed spaces equivalent of
\eq{ideq}. We have:

% Now, Lipschitz continuity
%implies, via Rademacher's theorem, that there exist sets of
%full-measure $\mcO_1\subset\pM_1$ and $\mcO_2\subset \pM_2$ on
%which $\phto$ and $\phot$ are differentiable. Set
%$$\mcU_1:=\mcO_1\cap \phot(\mcO_2)\;;$$
%Lipschitz continuity of $\phot$ implies that $\mcU_1$ also has
%full measure.

\begin{lemm}\label{Lnuisance} The map  $\phto$ is $C^1$
up-to-boundary.
\end{lemm}

\proof We can conformally rescale $g$ so that \eq{cm2h} holds with
$\tg=g_1$ and $\Omega=\Omega_1$; as \eq{cm2h} is conformally
invariant, \eq{cm2h} will also hold with $\tg=g_2$ and
$\Omega=\Omega_2$. Introducing coordinates $r_a$ as in the proof
of Theorem~\ref{Tmcc} with, say $k=0$, we can apply Theorem 3.3 of
\cite{ChNagyATMP} to obtain the desired conclusion. \qed

\medskip

Returning to the proof of Theorem~\ref{Tucc},
Lemma~\ref{Lnuisance} shows that for all $p\in \pM_1$ the maps
$(\phto){}_*(p)$ are similarities with non-zero ratio (in general
depending upon $p$). Differentiability of $\phto$ further implies
that $\phto$ is $ACL^n$, as defined in~\cite{Lelong-Ferrand}.   We
can then use a deep result of
Lelong-Ferrand~\cite[Theorem~A]{Lelong-Ferrand} to conclude that
$\varphi_2|_{\pM_1}$ and $\phto|_{\pM_1}$ are smooth. Now,
 $u:= (\varphi_2)^{\frac{n-2}{2}}$
solves the Yamabe equation, \bel{Yeq} \Delta_{g_1} u - \frac
{n-2}{4(n-1)} R_{g_1} u = (R_{g_2}\circ\phto)
u^{(n+2)/(n-2)}\;.\ee Here, as before, $R_{g_a}$ denotes the
curvature  scalar of the metric $g_a$. The right-hand-side of this
equation is in $L^\infty(\bM_1)$, and standard results on the
Dirichlet problem imply that $u$ --- and hence $\varphi_2$ --- is
uniformly $C^1$ on $\bM_1$. Now, $\phto$ is an isometry between
$g_1$ and $\varphi^{-2}_2g_2$ which implies that, in local
coordinates, $\phto$ satisfies on $M_1$  the over-determined set
of equations
$$\frac{\partial ^2\phto^i}{\partial x^\ell\partial x^m} =
\Gamma^k_{\ell m}(x) \frac{\partial \phto^i}{\partial x^k}-
\overline\Gamma^i_{rs}(\phto(x))\frac{\partial \phto^r}{\partial
x^\ell}\frac{\partial \phto^s}{\partial x^m}\;,
$$
where $\overline\Gamma^i_{rs}$ are, in local coordinates, the
Christoffel symbols of the Riemannian metric $(\varphi_2\circ\phot)^2 g_2$.
The right-hand-side of this set of equations extends by continuity
to a continuous function on $\bM_1$, which shows that $\phto$ is
uniformly $C^2$ on $\bM_1$. It follows that the right-hand-side of
\Eq{Yeq} is uniformly $C^1$ on $\bM_1$, hence $\varphi_2$ is
uniformly $C^2$ on $\bM_1$. An inductive repetition of this
argument establishes our claims. \qed

\begin{remk}\label{remk62}
It would be of interest to find a proof of Lemma~\ref{Lnuisance}
which is more in the spirit of conformal geometry than the methods
of~\cite{ChNagyATMP}.
\end{remk}

 \begin{remk} Recall that there exists a
conformally invariant version of the Abbott--Deser mass, due to
Ashtekar and Magnon~\cite{AshtekarMagnonAdS}, in a Lorentzian
space-time setting. We expect this expression to have a Riemannian
counterpart which is also conformally invariant, with the
numerical value thereof identical to that of  the expression we
propose. If that is the case, Theorem~\ref{Tmcc} is actually a
straightforward corollary of Theorem~\ref{Tucc}; in particular one
would not need to invoke the rather messy calculations of
Theorem~\ref{Tinv}, which are implicitly used in the proof
of Theorem~\ref{Tmcc}.
\end{remk}

\noindent{\bf{Acknowledgements.}} The authors are grateful to
Maung~Min-Oo and  Abdelghani~Zeghib  for useful discussions and to
the referee for his comments which led to an improved version of
Theorem \ref{Tmb}. PTC wishes to thank the Albert Einstein
Institute, Golm, for friendly hospitality during part of the work
on this paper.
%
%\bibliographystyle{amsplain}
%%\bibliographystyle{/usr/share/texmf/tex/revtex/prsty}
%\bibliography{%ptjjk,

\begin{thebibliography}{10}

\bibitem{AndDahl}
L.~Andersson and M.~Dahl, \emph{Scalar curvature rigidity for
asymptotically
  locally hyperbolic manifolds}, Annals of Global Anal.\ and Geom. \textbf{16}
  (1998), 1--27, dg-ga/9707017.

\bibitem{AshtekarMagnonAdS}
A.~Ashtekar and A.~Magnon, \emph{Asymptotically anti--de {S}itter
space-times},
  Class.\ Quantum Grav. \textbf{1} (1984), L39--L44.

\bibitem{BTZ}
M.~Ba{\~n}ados, C.~Teitelboim, and J.~Zanelli, \emph{{Black hole
in
  three-dimensional spacetime}}, Phys.\ Rev.\ Lett. \textbf{69} (1992),
  1849--1851.

\bibitem{Bartnik}
R.~Bartnik, \emph{The mass of an asymptotically flat manifold},
Commun.\ Pure and
  Appl.\ Math.\ \textbf{39} (1986), 661--693.

\bibitem{BartnikChrusciel}
R.~Bartnik and P.T. Chru\'sciel, \emph{Spectral boundary
conditions for
  {D}irac--type equations}, in preparation.

\bibitem{baum-killing}
H.~Baum, \emph{Complete {R}iemannian manifolds with imaginary
{K}illing
  spinors}, Ann. Glob. Anal. Geom \textbf{7} (1989), 205--226.

\bibitem{baum-dim5}
\bysame, \emph{Odd-dimensional {R}iemannian manifolds with
imaginary {K}illing
  spinors}, Ann.\ Glob.\ Anal.\ Geom. \textbf{7} (1989), 141--154.

\bibitem{Besse}
A.L. Besse, \emph{{Einstein manifolds}}, Ergebnisse der Mathematik
und ihrer
  Grenzgebiete. 3. Folge, vol.~10, Springer Verlag, Berlin, New York,
  Heidelberg, 1987.

\bibitem{BGH}
W.~Boucher, G.W. Gibbons, and G.T. Horowitz, \emph{Uniqueness
theorem for
  anti--de {S}itter spacetime}, Phys.\ Rev.\ \textbf{D30} (1984), 2447--2451.

\bibitem{Bourguignon92}
J-P. Bourguignon and P.~Gauduchon, \emph{Spineurs, op\'erateurs de
{D}irac et
  variations de m\'etriques}, Commun. Math. Phys. \textbf{144} (1992),
  581--599.

\bibitem{Bray:preparation2}
H.~Bray, \emph{Proof of the {R}iemannian {P}enrose conjecture
using the
  positive mass theorem}, Jour.\ Diff.\ Geom. \textbf{59} (2001), 177--267,
  math.DG/9911173.

\bibitem{BLP}
D.~Brill, J.~Louko, and P.~Peldan, \emph{Thermodynamics of
(3+1)-dimensional
  black holes with toroidal or higher genus horizons}, Phys.\ Rev. \textbf{D56}
  (1997), 3600--3610, gr-qc/9705012.

\bibitem{cheeger-LNM} J.~Cheeger, \emph{Critical points of distance functions and
applications to geometry}, Geometric topology: recent developments
(Montecatini Terme, 1990), Lect. Notes in Math \textbf{1504},
Springer, 1992.

\bibitem{ChBlesHouches}
Y.~Choquet-Bruhat, \emph{Positive-energy theorems}, Relativity,
groups and
  topology, II (Les Houches, 1983) (B.S. deWitt and R.~Stora, eds.),
  North-Holland, Amsterdam, 1984, pp.~739--785.

\bibitem{ChErice}
P.T. Chru\'sciel, \emph{Boundary conditions at spatial infinity
from a
  {H}amiltonian point of view}, Topological Properties and Global Structure of
  Space--Time (P.\ Bergmann and V.\ de~Sabbata, eds.), Plenum Press, New York,
  1986, pp. 49--59, {URL} \url{http://www.phys.univ-tours.fr/~piotr/scans}.

\bibitem{ChNagyATMP}
P.T. Chru\'sciel and G.~Nagy, \emph{The mass of spacelike
hypersurfaces in
  asymptotically {anti-de Sitter} space-times}, Adv.\ Theor.\ Math.\ Phys.
  \textbf{5} (2002), 697--754, gr-qc/0110014.

\bibitem{ChruscielSimon}
P.T. Chru\'sciel and W.~Simon, \emph{Towards the classification of
static
  vacuum spacetimes with negative cosmological constant}, Jour.\ Math.\ Phys.
  \textbf{42} (2001), 1779--1817, gr-qc/0004032.

\bibitem{Delay}
E.~Delay, \emph{Analyse pr\'ecis\'ee d'\'equations
semi-lin\'eaires elliptiques
  sur l'espace hyperbolique et application \`a la courbure scalaire conforme},
  Bull. Soc. Math. France \textbf{125} (1997), 345--381.

\bibitem{frkath}
T.~Friedrich and I.~Kath, \emph{Einstein manifolds of dimension
five with small
  first eigenvalue of the {D}irac operator}, Jour.\ Diff.\ Geom. \textbf{29}
  (1989), 263--279.

\bibitem{pg-italien}
P.~Gauduchon, \emph{Hermitian connections and {D}irac operators},
Bolletino U.
  M. I. \textbf{11 B} (1997), no.~suppl., 257--288.

\bibitem{GHHP}
G.W. Gibbons, S.W. Hawking, G.T. Horowitz, and M.J. Perry,
\emph{Positive mass
  theorem for black holes}, Commun. Math. Phys. \textbf{88} (1983), 295--308.

\bibitem{mh-inegalite-penrose}
M.~Herzlich, \emph{A {Penrose-}like inequality for the mass of
{R}iemannian
  asymptotically flat manifolds}, Commun. Math. Phys. \textbf{188} (1997),
  121--133.

\bibitem{Herzlich:mass}
\bysame, \emph{The positive mass theorem for black holes
revisited}, Jour.\
  Geom.\ Phys. \textbf{26} (1998), 97--111.

\bibitem{HMZ-asian}
O.~Hijazi, S.~Montiel, and X.~Zhang, \emph{Conformal lower bounds
for the Dirac operator of embedded hypersurfaces}, Asian J. Math.
\textbf{6} (2002), 23--36.

\bibitem{HorowitzMyers}
G.T. Horowitz and R.C. Myers, \emph{The {AdS/CFT} correspondence
and a new
  positive energy conjecture for general relativity}, Phys. Rev. \textbf{D59}
  (1999), 026005 (12 pp.).

\bibitem{HI1}
G.~Huisken and T.~Ilmanen, \emph{The {Riemannian P}enrose
inequality}, Int.\
  Math.\ Res.\ Not. \textbf{20} (1997), 1045--1058.

\bibitem{HI2}
\bysame, \emph{The inverse mean curvature flow and the {Riemannian
P}enrose
  inequality}, Jour.\ Diff.\ Geom. \textbf{59} (2001), 353--437, {URL}
  \url{http://www.math.nwu.edu/~ilmanen}.

\bibitem{Kottler}
F.~Kottler, \emph{{\"Uber die physikalischen Grundlagen der
Einsteinschen
  Gravitationstheorie}}, Annalen der Physik \textbf{56} (1918), 401--462.

\bibitem{Kuehnel}
W.~Kuehnel, \emph{{Conformal transformations between Einstein
spaces}},
  Conformal geometry (R.S. Kulkarni and U.~Pinkall, eds.), Aspects Math.: E,
  12, F.~Vieweg \& Sohn, Braunschweig, 1988, pp.~105--146.

\bibitem{lebrun-cex}
C.~{Le Brun}, \emph{Counterexamples to the generalized positive
action
  conjecture}, Commun. Math. Phys. \textbf{118} (1988), 591--596.

\bibitem{Lelong-Ferrand}
J.~Lelong-Ferrand, \emph{{Geometrical interpretations of scalar
curvature and
  regularity of conformal homeomorphisms }}, {Differ. Geom. Relativ., Vol.
  Honour A. Lichnerowicz 60th Birthday}, 1976, pp.~91--105.

\bibitem{LiTam}
P.~Li and L.-F. Tam, \emph{Complete surfaces with finite total
curvature},
  Jour.\ Diff.\ Geom. \textbf{33} (1991), 139--168.

\bibitem{Lichnerowicz63}
A.~Lichnerowicz, \emph{Spineurs harmoniques}, C.R. Acad. Sci.
Paris S\'er. A-B
  \textbf{257} (1963), 7--9.

\bibitem{min-oo:scalar}
M.~Min-Oo, \emph{Scalar curvature rigidity of asymptotically
hyperbolic spin
  manifolds}, Math.\ Ann. \textbf{285} (1989), 527--539.

\bibitem{ParkerTaubes82}
T.~Parker and C.~Taubes, \emph{On {W}itten's proof of the positive
energy
  theorem}, Commun. Math. Phys. \textbf{84} (1982), 223--238.

\bibitem{Schrodinger32}
E.~Schr\"odinger, \emph{Diracsches {E}lektron im {S}chwerfeld},
Preuss. Akad.
  Wiss. Phys.-Math. \textbf{11} (1932), 436--460.

\bibitem{Shiohama}
K.~Shiohama, \emph{Total curvature and minimal area of complete
open surfaces},
  Proc.\ Am.\ Math.\ Soc. \textbf{94} (1985), 310--316.

\bibitem{Vanzo}
L.~Vanzo, \emph{Black holes with unusual topology}, Phys. Rev.
\textbf{D56}
  (1997), 6475--6483, gr-qc/9705004.

\bibitem{Wang}
X.~Wang, \emph{Mass for asymptotically hyperbolic manifolds},
Jour.\ Diff.\
  Geom. \textbf{57} (2001), 273--299.

\bibitem{witten:mass}
E.~Witten, \emph{A simple proof of the positive energy theorem},
Commun. Math.
  Phys. \textbf{80} (1981), 381--402.

\bibitem{Zhang:hpet}
X.~Zhang, \emph{A definition of total energy-momenta and the
positive mass
  theorem on asymptotically hyperbolic 3 manifolds {I}},  (2001), preprint.

\end{thebibliography}
%../../references/newbiblio,%
%../../references/reffile,%
%../../references/bibl,%
%../../references/Energy,%
%../../references/hip_bib,%
%../../references/netbiblio}
%%\texttt{\input{READMEl}}
\def\cprime{$'$}
\providecommand{\bysame}{\leavevmode\hbox
to3em{\hrulefill}\thinspace}
\providecommand{\MR}{\relax\ifhmode\unskip\space\fi MR }
% \MRhref is called by the amsart/book/proc definition of \MR.
\providecommand{\MRhref}[2]{%
  \href{http://www.ams.org/mathscinet-getitem?mr=#1}{#2}
} \providecommand{\href}[2]{#2}

\end{document}